\documentclass[11pt]{amsart}
\usepackage{amsmath,amssymb, amsbsy}
\usepackage{color,psfrag}
\usepackage[dvips]{graphicx}

\usepackage{enumerate}
\newcommand{\R}{{\mathbb R}}
\newcommand{\N}{{\mathbb N}}
\newcommand{\Z}{{\mathbb Z}}

\newcommand{\eps}{\varepsilon}

\newcommand{\E}{{\mathcal E }}
\renewcommand{\phi}{\varphi}
\def\endproof{\hfill $\Box$\par\vskip3mm}

\def\ds{\displaystyle}
\renewcommand{\ge }{\geqslant}

\renewcommand{\le }{\leqslant}
\renewcommand{\leq }{\leqslant}
\def\neweq#1{\begin{equation}\label{#1}}
\def\endeq{\end{equation}}
\def\eq#1{(\ref{#1})}

\newcommand{\pt}[1]{\left( #1 \right) }

\newtheorem{theorem}{Theorem}
\newtheorem{proposition}[theorem]{Proposition}
\newtheorem{lemma}[theorem]{Lemma}
\newtheorem{corollary}[theorem]{Corollary}
\newtheorem{conjecture}[theorem]{Conjecture}

\newtheorem{remark}[theorem]{Remark}
\newtheorem{definition}{Definition}

\begin{document}

\title{Energy transfer between modes in a nonlinear beam equation}


\author[Ubertino BATTISTI]{Ubertino BATTISTI}
\address{\hbox{\parbox{5.7in}{\medskip\noindent{Dipartimento di Scienze Matematiche, \\
Politecnico di Torino,\\
        Corso Duca degli Abruzzi 24, 10129 Torino, Italy. \\[3pt]
        \em{E-mail address: }{\tt ubertino.battisti@polito.it}}}}}
\author[Elvise BERCHIO]{Elvise BERCHIO}
\address{\hbox{\parbox{5.7in}{\medskip\noindent{Dipartimento di Scienze Matematiche, \\
Politecnico di Torino,\\
        Corso Duca degli Abruzzi 24, 10129 Torino, Italy. \\[3pt]
        \em{E-mail address: }{\tt elvise.berchio@polito.it}}}}}
\author[Alberto FERRERO]{Alberto FERRERO}
\address{\hbox{\parbox{5.7in}{\medskip\noindent{Dipartimento di Scienze e Innovazione Tecnologica, \\
Universit\`a del Piemonte Orientale ``Amedeo Avogadro'',\\
        Viale Teresa Michel 11, 15121 Alessandria, Italy. \\[3pt]
        \em{E-mail address: }{\tt alberto.ferrero@uniupo.it}}}}}
\author[Filippo GAZZOLA]{Filippo GAZZOLA}
\address{\hbox{\parbox{5.7in}{\medskip\noindent{Dipartimento di Matematica,\\
Politecnico di Milano,\\
   Piazza Leonardo da Vinci 32, 20133 Milano, Italy. \\[3pt]
        \em{E-mail address: }{\tt filippo.gazzola@polimi.it}}}}}
\maketitle

\begin{abstract}
We consider the nonlinear nonlocal beam evolution equation introduced by Woinowsky-Krieger \cite{woinowsky}. We study the existence and behavior of periodic solutions: these are called nonlinear modes. Some solutions only have two active modes and we investigate
whether there is an energy transfer between them. The answer depends on the geometry of the energy function which, in turn, depends on the amount of
compression compared to the spatial frequencies of the involved modes. Our results are complemented with numerical experiments; overall, they
give a complete picture of the instabilities that may occur in the beam. We expect these results to hold also in more complicated dynamical
systems.\par\noindent
{\bf R\'esum\'e:} On consid\`ere l'\'equation d'\'evolution de la poutre nonlin\'eaire et nonlocale introduite par Woinowsky-Krieger \cite{woinowsky}. On \'etudie
l'existence et le comportement des solutions p\'eriodiques: on les appelle modes nonlin\'eaires. Certaines solutions ont seulement deux modes actifs et nous
\'etudions le possible transfer d'\'energie entre eux. La r\'eponse d\'epend de la g\'eom\'etrie de la fonctionnelle d'\'energie qui, \`a son tour, d\'epend
de la quantit\'e de compression et des fr\'equences spatiales des modes actifs. Nos r\'esultats sont compl\'et\'es par des experiments num\'eriques; ils donnent
une description d'ensemble assez compl\`ete des instabilit\'es qui peuvent appara\^itre dans la poutre. On s'attend \`a ce que ces r\'esultats soient valables aussi
pour des syst\`emes dynamiques plus compliq\'es.\par\noindent
\end{abstract}

\section{Introduction}

In 1950, Woinowsky-Krieger \cite{woinowsky} modified the classical beam models by Bernoulli and Euler assuming a nonlinear dependence of the axial
strain on the deformation gradient, by taking into account the stretching of the beam due to its elongation. Let us mention that, independently,
Burgreen \cite{burg} derived the very same nonlinear beam equation which reads
$$
M\, u_{tt}+EI\, u_{xxxx}+\Big[P-\eta\, \|u_x\|^2_{L^2(0,\ell)}\Big]u_{xx}=f\qquad x\in(0,\ell)\, ,\ t>0\, ,
$$
where $u$ denotes the vertical displacement of the beam whose length is $\ell$. The constant $\eta>0$ depends on the elasticity of the material composing
the beam and the term $\eta\|u_x\|^2_{L^2(0,\ell)}$ measures the geometric nonlinearity of the beam due to
its stretching. The constant $P$ is the axial force acting at the endpoints of the beam: a positive $P$ means that the beam is compressed while
a negative $P$ means that the beam is stretched. We are mainly interested in compressed beams ($P>0$) although some of our results also apply to
free ($P=0$) and stretched ($P<0$) beams. Finally, $M>0$ denotes the mass per unit length, $EI>0$ is the flexural rigidity of the beam,
whereas $f=f(x,t)$ is an external load.\par
We assume that the beam is hinged at its endpoints and this results in the so-called Navier boundary conditions.
For simplicity, we consider a beam lying on the segment $x\in(0,\pi)$, we normalize the constants, we take null force, and we reduce to
\neweq{truebeam}
\left\{\begin{array}{ll}
u_{tt}+u_{xxxx}+\Big[P-\frac{2}{\pi}\, \|u_x\|^2_{L^2(0,\pi)}\Big]u_{xx}=0\quad & x\in(0,\pi)\, ,\ t>0\, ,\\
u(0,t)=u(\pi,t)=u_{xx}(0,t)=u_{xx}(\pi,t)=0\quad & t>0\, .
\end{array}\right.
\endeq
A description of \eqref{truebeam} with $f\neq 0$ would require a huge effort and falls beyond the scopes of this paper. We expect this kind of analysis to require the exploitation of previous results for related forced ODE's, see for instance \cite{burza,dingza}. The existence and uniqueness of global solutions of the initial value problem associated to \eq{truebeam} has been proved in \cite{ball,dickey0}, while in \cite{grotta,yaga1} the existence of chaotic dynamics for  \eqref{truebeam} was shown. In this paper we perform a detailed (theoretical and numerical) study of the stability of its nonlinear modes. We make use of refined properties of the Hill and Duffing equations, classical tools from Floquet theory (such as the monodromy matrices and Poincar\'e maps), some stability criteria and estimates of elliptic functions, and numerical experiments when these theoretical arguments fail. Let us describe our results.

\begin{figure}[ht]
\begin{center}
{\includegraphics[height=35mm, width=65mm]{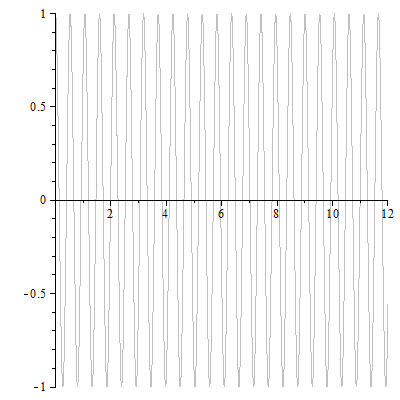}}\qquad{\includegraphics[height=35mm, width=65mm]{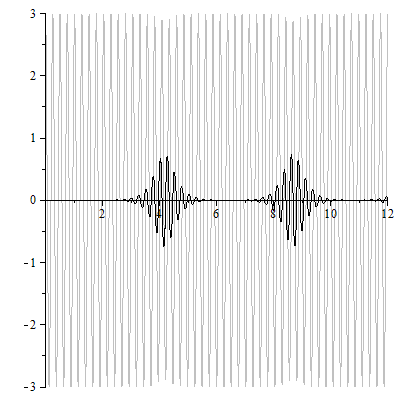}}
\caption{Stable (left) and unstable (right) oscillations.}\label{typical}
\end{center}
\end{figure}

It is well-known that the unforced evolution equation \eq{truebeam} admits infinitely many {\em nonlinear modes}, that is, solutions having a unique nontrivial periodic-in-time Fourier component, see Definition \ref{nonlinmod}.
The Fourier component is the solution of a Duffing equation \cite{duffing} which is obtained from \eq{truebeam} by separating variables. The behavior of the Duffing equation changes if the compression parameter $P$ is above or below a threshold which depends on the considered Fourier component. In Theorems \ref{periodicT1} and \ref{periodicT2} we analyze with great precision the dependence of the period and of the amplitude of the solutions of the Duffing equations with respect to the internal energy of the beam.
Then we enter into the main core of the paper. Local nonlinear wave equations admit infinitely many resonances
since the dynamical system itself is infinite dimensional, see \cite[$\S E.3.4$]{SVM}: indeed, in these equations the initial energy of the system immediately spreads on infinitely many modes and therefore the resonances are difficult to detect; see e.g.\ \cite{bfg1} for a plate equation. On the contrary, for nonlocal equations such
as \eq{truebeam} the energy may remain confined to a finite number of modes.
Recent results in \cite{bfgk,gazkar,gazpav,gazpav2,gazpav3} highlight unexpected amplifying
oscillations in stationary nonlinear beam equations and, in this paper, we aim to study whether these oscillations transfer from one mode to another.
We consider particular solutions of \eq{truebeam} which only have two nontrivial time-dependent Fourier coefficients, one being initially smaller than the other by several orders of magnitude. We study the stability of the large mode with respect to the small mode. The typical pictures describing the loss of stability are as in Figure \ref{typical}.
In both pictures, the gray oscillations represent the large mode whereas the black oscillations represent the small mode. In the left picture, the initial data are such that no black oscillations are visible, which means that the large mode is stable. In the right picture, we increase the initial data and one
may see a large oscillation also in the small mode: this mode suddenly grows up by capturing some energy from the large mode which
decreases its amplitude of oscillation when the transfer of energy occurs. This is what we call instability of a mode with respect to another mode: the instability manifests through a sudden transfer of energy between modes. Since the frequency of a nonlinear mode
depends on the amplitude of oscillation (and hence on the energy), in some cases the energy transfer may or may not occur according to the amount of energy inside
the system.
\par
As pointed out by Stoker in \cite[Chapter IV]{stoker}, \emph{in any consideration of stability of a given system one fundamental difficulty is that of defining the notion of stability in a logical and reasonable manner without destroying the chances of applying the definition in a practical way.}
In this paper we deal with the linear stability which is characterized in Definition \ref{defstabb}.

Not only the stability analysis depends on the modes considered, but it also strongly depends on the magnitude of the compression $P$
and several different cases have to be analyzed, according to the value of $P$ with respect to the spatial frequencies of the two modes involved.
In some situations we take advantage of the stability study for large energies due to Cazenave-Weissler \cite{cazw,cazw2}, in some other cases
we use some stability criteria (recalled in Section \ref{stabcrit}) for the Hill equation, further cases require ``by hand'' estimates.
We complement the theoretical results with numerical experiments; this leads to conjectures and to several open problems.
Each proof has its own difficulties but two of them are particularly involved, those of Theorems \ref{stability11} and \ref{stability22}.
The former combines stability criteria for the Hill equation with delicate properties of elliptic integrals, whereas the latter makes use of
fine asymptotic estimates for the solutions of a Duffing equation with negative energies.\par
A further motivation for this paper is to give some hints about the nonlinear structural behavior of suspension bridges
\cite{bfg1,gazz,bookgaz}: it is reasonable to expect that if some instability appears in a simplified model such as \eq{truebeam}, namely if the deck of the bridge is seen as a beam, then similar instabilities will appear in more sophisticated models. One may take advantage of the explicit solutions and of the precise results that we reach for \eq{truebeam} in order to guess some responses for more complicated dynamical systems. Our results clearly show that the transfer of energy between modes depends on the ratio of their spatial frequencies.
Some couples of modes never transfer energy to each other while some different couples are more prone to an energy transfer.
The energy threshold of instability depends on the considered couple and if one aims to prevent some particular dangerous oscillations, such as torsional oscillations in plates modeling suspension bridges, one should also prevent the appearance of those oscillations which are prone to transfer their
energy to the dangerous ones.

\section{Nonlinear modes}\label{nonlinfor}
For the main properties of the stationary solutions to \eq{truebeam}, see Proposition \ref{homogeneous} in the Appendix. We discuss here some basic facts related to the evolution equation \eq{truebeam}. We refer to \cite{burg,eisley,stoker} for former works on this topic. We state and prove all the results in detail because we need very precise statements for the stability analysis in the subsequent sections.

We first characterize the nonlinear modes of \eq{truebeam}
by considering solutions in the form
\neweq{form}
v_k(x,t)=\Theta_k(t)\sin(kx)\,.
\endeq
\begin{definition}\label{nonlinmod}
We call a function $v_k$ in the form \eqref{form} a $k$-th nonlinear mode of \eqref{truebeam}.
\end{definition}
It is straightforward that $v_k$ in \eq{form} satisfies the boundary conditions in \eq{truebeam}. Furthermore, by inserting \eq{form} into \eq{truebeam}, it is readily seen that the Fourier coefficient $\Theta_k$ satisfies
\neweq{ODE}
\ddot{\Theta}_k(t)+k^2(k^2-P)\Theta_k(t)+k^4\Theta_k(t)^3=0\qquad(t>0)
\endeq
and its behavior depends on whether $k^2\lessgtr P$.
When $k^2-P>0$, \eq{ODE} is the so-called Duffing equation which was introduced in \cite{duffing} to describe a nonlinear oscillator with a cubic stiffness, see also \cite{stoker}. The name Duffing equation is nowadays also attributed to \eq{ODE} when the coefficient of the linear term is nonpositive, see \cite[Section 2.2]{guck}. To \eq{ODE} we associate some initial values
\neweq{alphabeta}
\Theta_k(0)=\alpha\ ,\qquad \dot{\Theta}_k(0)=\beta\, ,\qquad(\alpha,\beta\in\R)\,
\endeq
and the corresponding constant energy:
\neweq{energyy}
E(\alpha,\beta)=\frac{\dot{\Theta}_k^2}{2}+\frac{k^2(k^2-P)}{2}\Theta_k^2+\frac{k^4}{4}\Theta_k^4\equiv\frac{\beta^2}{2}+\frac{k^2(k^2-P)}{2}\alpha^2+
\frac{k^4}{4}\alpha^4\,.
\endeq
For all $E>0$ we put
\neweq{Lambdas}
\Lambda_1(E):=\frac{\sqrt{(k^2-P)^2+4E}+P-k^2}{k^2}>0\, ,\quad\Lambda_2(E):=\frac{\sqrt{(k^2-P)^2+4E}-P+k^2}{k^2}>0\, ,
\endeq
while for $k^2<P$ and $-\frac{(P-k^2)^2}{4}<E<0$ we define
\neweq{delta}
\delta:=\frac{P-k^2-\sqrt{(P-k^2)^2+4E}}{2\sqrt{|E|}}\, .
\endeq
We point out that there exist infinitely many $k$-th nonlinear modes for each $k$ and that they are not proportional to each other.
Their shape is described by the solution $\Theta_k$ of \eq{ODE} which depends on the initial energy $E(\alpha,\beta)$ in \eq{energyy}.
For this reason, with an abuse of language, we will also call $\Theta_k$ a nonlinear mode of \eq{truebeam}. Some properties of the $\Theta_k$ that will be useful in the sequel are collected in Theorems \ref{periodicT1} and \ref{periodicT2} below. These results adapt to our context previous statements by Burgreen \cite{burg}. Since we also need some tools from
their proofs, we briefly sketch them in Section \ref{pT1}. The first statement deals with the beam under small compression.

\begin{theorem}\label{periodicT1}
Take an integer $k\ge1$ and let $P\in\R$ be such that $P\le k^2$. If a solution $\Theta_k$ of \eqref{ODE} takes the initial values \eqref{alphabeta} for
some $(\alpha,\beta)\neq (0,0)$, thereby satisfying
\neweq{lbforE}
E(\alpha,\beta)>0\, ,
\endeq
then $\Theta_k$ is periodic and its period is given by
\neweq{TE}
T(E)=\frac{4\sqrt2}{k^2}\int_0^1\frac{d\theta}{\sqrt{(\Lambda_2(E)+\Lambda_1(E)\theta^2)(1-\theta^2)}}\, ,
\endeq
see \eqref{Lambdas}. In particular, the map $E\mapsto T(E)$ is strictly decreasing on $(0,+\infty)$ and
\neweq{limitperiod}
\lim_{E\downarrow0}T(E)=\frac{2\pi}{k\sqrt{k^2-P}}\, (=+\infty\mbox{ if }k^2=P)\, ,\qquad\lim_{E\to+\infty}T(E)=0\, .
\endeq
\end{theorem}

The second statement deals with the beam under large compression.

\begin{theorem}\label{periodicT2}
Take an integer $k\ge1$ and let $P\in\R$ be such that $P>k^2$.
If a solution $\Theta_k$ of \eqref{ODE} takes the initial values \eqref{alphabeta} for some $\alpha,\beta\in\R$ satisfying
\neweq{lbforE2}
\mbox{either}\qquad-\frac{(P-k^2)^2}{4}<E(\alpha,\beta)<0\qquad\mbox{or}\qquad E(\alpha,\beta)>0\, ,
\endeq
then $\Theta_k$ is periodic and its period is given by \eqref{TE} if $E>0$ and by
\neweq{TE2}
T(E)=\frac{2\sqrt{2}}{k\sqrt{P-k^2+\sqrt{(P-k^2)^2+4E}}}\, \int_\delta^1\frac{d\theta}{\sqrt{(1-\theta^2)(\theta^2-\delta^2)}}
\endeq
if $E<0$, where $\delta\in(0,1)$ is defined in \eqref{delta}. In particular,
\neweq{limitperiod2}
\lim_{E\downarrow-\frac{(P-k^2)^2}{4}}T(E)=\frac{\pi\, \sqrt2}{k\, \sqrt{P-k^2}}\, ,\qquad
\lim_{E\to0}T(E)=+\infty\, ,\qquad\lim_{E\to+\infty}T(E)=0\, .
\endeq
Moreover, the map $E\mapsto T(E)$ is strictly increasing on $(-\frac{(P-k^2)^2}{4},0)$ and strictly decreasing on $(0,+\infty)$.
\end{theorem}

The period $T(E)$ can also be expressed differently from \eq{TE}, see \eq{newTE}. In Section \ref{some} of the Appendix we comment the
assumptions \eq{lbforE} and \eq{lbforE2}.\par\smallskip

\section{Stability of the nonlinear modes}

\subsection{Linear stability}

We consider here solutions of \eq{truebeam} having only two nontrivial Fourier components, that is:
\neweq{form2}
u(x,t)=w(t)\sin(mx)+z(t)\sin(nx)
\endeq
for some integers $n,m\ge1$, $n\neq m$. One does not expect such $u$ to be periodic-in-time but we will show that it may have both the tendency to become
periodic and to break down periodicity. After inserting \eq{form2} into \eq{truebeam} we reach the following (nonlinear) system:
\neweq{cw}
\left\{\begin{array}{l}
\ddot{w}(t)+m^2(m^2-P)w(t)+m^2\big(m^2w(t)^2+n^2z(t)^2\big)w(t)=0\, ,\\
\ddot{z}(t)+n^2(n^2-P)z(t)+n^2\big(m^2w(t)^2+n^2z(t)^2\big)z(t)=0\ ,
\end{array}\right.
\endeq
to which we associate the initial conditions
\neweq{initialsyst}
w(0)=w_0\, ,\ \dot{w}(0)=w_1\, ,\quad z(0)=z_0\, ,\ \dot{z}(0)=z_1\, .
\endeq
Also system \eq{cw} is conservative and its constant energy is
\begin{eqnarray*}
\E(w_0,w_1,z_0,z_1)\!&=&\!\tfrac{\dot{w}^2}{2}\!+\!\tfrac{\dot{z}^2}{2}\!+\!m^2(m^2\!-\!P)\tfrac{w^2}{2}\!+\!n^2(n^2\!-\!P)\tfrac{z^2}{2}\!+\!
m^4\tfrac{w^4}{4}\!+\!n^4\tfrac{z^4}{4}\!+\!m^2n^2\tfrac{w^2z^2}{2} \notag \\
\!&\equiv&\!\tfrac{w_1^2}{2}\!+\!\tfrac{z_1^2}{2}\!+\!m^2(m^2\!-\!P)\tfrac{w_0^2}{2}\!+\!n^2(n^2\!-\!P)\tfrac{z_0^2}{2}\!+\!
m^4\tfrac{w_0^4}{4}\!+\!n^4\tfrac{z_0^4}{4}\!+\!m^2n^2\tfrac{w_0^2z_0^2}{2}\, .
\end{eqnarray*}
This energy consists of three terms: the total energy $E_w$ of $w$ (kinetic+potential energy), the total energy $E_z$ of $z$, and the coupling energy $E_{wz}$.
Although their sum is constant, these three energies depend on time and they are explicitly given by
\neweq{energies3}
E_w(t)\!:=\!\tfrac{1}{2}\dot{w}^2\!+\!\tfrac{m^2(m^2-P)}{2}w^2\!+\!\tfrac{m^4}{4}w^4\, ,\quad
E_z(t)\!:=\!\tfrac{1}{2}\dot{z}^2\!+\!\tfrac{n^2(n^2-P)}{2}z^2\!+\!\tfrac{n^4}{4}z^4\, ,\quad
E_{wz}(t)\!:=\!\tfrac{m^2n^2}{2}w^2z^2\, .
\endeq

\begin{remark} {\em With the change of unknowns $w\mapsto mw$ and $z\mapsto nz$ the system \eq{cw} simplifies and reads
\neweq{cwx}
\left\{\begin{array}{l}
\ddot{w}(t)+m^2(m^2-P)w(t)+m^2\big(w(t)^2+z(t)^2\big)w(t)=0\, ,\\
\ddot{z}(t)+n^2(n^2-P)z(t)+n^2\big(w(t)^2+z(t)^2\big)z(t)=0\ .
\end{array}\right.
\endeq
But, in order to maintain the notations used so far, we remain with \eq{cw}. The system \eq{cwx} will be used in the proof of
Theorem \ref{stability121}.\endproof}\end{remark}

We wish to analyze the stability of the modes $\Theta_m$ (the solution of \eq{ODE}-\eq{alphabeta}) within the nonlinear system \eq{cw} in both
cases $m\lessgtr n$ and for $P\ge0$ being in different positions with respect to $m^2$ and $n^2$. The stability properties of $\Theta_m$ depend on the
energy of $\Theta_m$ which, in the sequel, will be denoted by $E_{\Theta_m}$. From \eq{energyy} we recall that
$$
E_{\Theta_m}=\frac{\dot{\Theta}_m(0)^2}2+\frac{m^2(m^2-P)}2 \Theta_m(0)^2+\frac{m^4}4 \Theta_m(0)^4 \, .
$$

\begin{definition}\label{defstabb}
The mode $\Theta_m$ is said to be {\bf linearly stable} ({\bf unstable}) with respect to the $n$-th mode $\Theta_n$ if $\xi\equiv0$ is a stable (unstable) solution of
the Hill equation
\neweq{hill2}
\ddot{\xi}+a(t)\xi=0\, ,\qquad a(t)=n^2(n^2-P)+m^2n^2\Theta_m(t)^2\quad \forall t\, .
\endeq
\end{definition}

There exist also stronger definitions of stability. The mode $\Theta_m$ is said to be {\bf orbitally stable} if for any $\eps>0$ there exists $\delta>0$ such that if $(w(t),z(t))$ is a solution of \eqref{cw}-\eqref{initialsyst} with
$$
\min_{s\in[0,T(E)]}\{|w_0-\Theta_m(s)|,|w_1-\dot{\Theta}_m(s)|,|z_0|,|z_1|\}<\delta\, ,
$$
then
$$
{\ds \sup_{t\in \R} \min_{s\in[0,T(E)]}
\big(\left|w(t)-\Theta_m(s)\right|+\big|\dot{w}(t)-\dot{\Theta}_m(s)\big|+|z(t)|+|\dot{z}(t)|\big)<\eps} \, .
$$

The mode $\Theta_m$ is said to be {\bf orbitally unstable} if it is not orbitally stable. In general, it is not true that linear stability implies orbital stability. In some cases, the two concepts are equivalent, see for example \cite[Theorems 2.5-2.6]{ghg} where, by exploiting the KAM theory, sufficient conditions for the equivalence of the two notions are provided. For system \eqref{cw} we prove that linear instability implies orbital instability, see the end of Section \ref{s:cw}. Moreover, a result by Ortega \cite{ortega} states that if the trivial solution $\xi\equiv0$ of \eq{hill2} is stable, then also the trivial solution
of the {\em nonlinear} Hill equation
$$\ddot{\xi}+a(t)\xi+n^4\xi^3=0\, ,\qquad a(t)=n^2(n^2-P)+n^2m^2 w(t)^2\quad \forall t$$
is stable. Therefore, the linear stability appears to be a satisfactory definition also in nonlinear regimes.

Since \eq{hill2} is linear, the linear stability is equivalent to state that all the solutions of \eq{hill2} are bounded. On the other hand, since \eq{ODE}
is nonlinear, the stability of $\Theta_m$ depends on the initial conditions \eq{alphabeta} and on the corresponding energy \eq{energyy}.
On the contrary, the linear instability of $\Theta_m$ occurs when the trivial solution of \eq{hill2} is unstable. This means that if we consider
a solution of \eq{cw} with $|z(0)|+|\dot{z}(0)|\ll|w(0)|+|\dot{w}(0)|$, that is, the initial energy is almost all due to the term $E_w(0)$, the component
$w$ conveys part of its energy to $z$ for $t>0$, see Figure \ref{typical}.
\par

The potential energy of the system \eq{cw} is given by
\neweq{potential}
{\mathcal U }(w,z)=m^2(m^2-P)\frac{w^2}{2}+n^2(n^2-P)\frac{z^2}{2}+m^4\frac{w^4}{4}+n^4\frac{z^4}{4}+m^2n^2\frac{w^2z^2}{2}\, .
\endeq
The orbits of \eq{cw} lie inside the sublevels of ${\mathcal U }$; if $\E_0$ denotes the initial (and constant) energy of \eq{cw}, one has
${\mathcal U }(w(t),z(t))\le\E_0$ for all $t\ge0$. The function ${\mathcal U }$ has different geometries according to the mutual positions of $m^2$, $n^2$, and $P$.
Below we state our stability results by discussing separately the three different geometries of ${\mathcal U }$. For the cases not covered by our theoretical statements, we numerically compute the
eigenvalues of the Poincar\'e map of the linearized system. See Section \ref{s:cw} for the relation between these eigenvalues and linear stability. The numerical results suggest a number of conjectures that we put near to the corresponding theoretical statements.

\subsection{The convex case}

\begin{figure}[ht]
\begin{center}
 {\includegraphics[height=44mm, width=76mm]{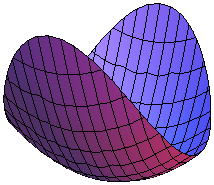}}
\caption{The potential energy functional ${\mathcal U }$ when $0\le
P\le\min\{m^2,n^2\}$.}\label{Ppiccolo}
\end{center}
\end{figure}

We consider first the case where $0\le P\le \min\{m^2,n^2\}$. Then the functional ${\mathcal U }$ in \eq{potential} is convex and its qualitative graph is
plotted in Figure \ref{Ppiccolo}. All the sublevels of ${\mathcal U }$ resemble to ellipses.

In this case, the largest mode is stable for both small and large energies.

\begin{theorem}\label{stability0}
Assume that $0\le P\le n^2<m^2$. Then there exist $0<E_1\le E_2<\infty$ such that
$\Theta_m$ is linearly stable with respect to $\Theta_n$ whenever $0<E_{\Theta_m}\le E_1$ or $E_{\Theta_m}>E_2$.
\end{theorem}
Numerical results suggest the following

\begin{conjecture}\label{ovvia}
Assume that $0\le P\le n^2<m^2$. Then $\Theta_m$ is linearly stable with respect to $\Theta_n$ for all $E_{\Theta_m}>0$.
\end{conjecture}

In favor of this conjecture we also have the following statement.

\begin{theorem}\label{stability11}
Assume that $0\le P\le n^2\le m^2(\tfrac{21}{22})^2$. Then $\Theta_m$ is linearly stable with respect to $\Theta_n$ for all $E_{\Theta_m}>0$.
\end{theorem}

The proof of Theorem \ref{stability11} requires a careful analysis of the complete elliptic integral of the first kind $K(\cdot)$
and its comparison with some functions arising from the stability criteria in Proposition \ref{lyapzhu} in the Appendix.
Theorem \ref{stability11} holds under the assumption that $m\ge1.04762n$, which is stronger than the optimal assumption $m>n$.
From Theorem \ref{stability11} we infer that Conjecture \ref{ovvia} is true for small modes.

\begin{corollary}\label{conj}
If $m\le22$ and $0\le P\le n^2<m^2$, then $\Theta_m$ is linearly stable with respect to $\Theta_n$ for all $E_{\Theta_m}>0$.
\end{corollary}

In fact, our proofs may be extended to $m\le27$, see Remark \ref{better} at the end of the proof of Theorem \ref{stability11};
in other words, Conjecture \ref{ovvia} holds under the additional assumption that $m\ge1.03847n$.\par
The stability analysis is fairly different as far as the mode with smaller frequency is involved. We define the two sets
\begin{align}\label{eq:condinst}
I_U:=\bigcup_{k\in\N}\Big(\pt{k+1}\pt{2k+1},\pt{k+1}\pt{2k+3}\Big)\, ,\quad I_S:=\bigcup_{k\in\N}\Big(k\pt{2k+1},\pt{k+1}\pt{2k+1}\Big)\, .
\end{align}
Note that $\overline{I_S\cup I_U}=[0,+\infty)$. Then we prove

\begin{theorem}\label{stability12}
Assume that $0\le P<m^2<n^2$, let $I_U$ and $I_S$ be as in \eqref{eq:condinst}. There exist $0<E_1\le E_2$ such that:
\begin{itemize}
\item[(i)] if $0<E_{\Theta_m}<E_1$ then $\Theta_m$ is linearly stable with respect to $\Theta_n$;
\item[(ii)] if $\frac{n^2}{m^2}\in I_U$ and $E_{\Theta_m}>E_2$ then $\Theta_m$ is linearly unstable with respect to $\Theta_n$;
\item[(iii)] if $\frac{n^2}{m^2}\in I_S$ and $E_{\Theta_m}>E_2$ then $\Theta_m$ is linearly stable with respect to $\Theta_n$.

\end{itemize}
\end{theorem}

In the limit case $0<P=m^2<n^2$, the following result for large energies holds
\begin{theorem} \label{t:limit-case}
Assume that $0<P=m^2<n^2$, let $I_U$ and $I_S$ be as in \eqref{eq:condinst}. There exists $E_2>0$ such that: \par
\begin{itemize}
\item[(i)] if $\frac{n^2}{m^2}\in I_U$ and $E_{\Theta_m}>E_2$ then $\Theta_m$ is linearly unstable with respect to $\Theta_n$;
\item[(ii)] if $\frac{n^2}{m^2}\in I_S$ and $E_{\Theta_m}>E_2$ then $\Theta_m$ is linearly stable with respect to $\Theta_n$.
\end{itemize}
\end{theorem}
Theorem \ref{t:limit-case} does not deal with stability for small energies. The numerical experiments
that we performed in this case suggest the following

\begin{conjecture}
Assume that $0<P=m^2<n^2$ and let $E_2$ be as in Theorem \ref{t:limit-case}. There exists $0<E_1\le E_2$ such that if $0<E_{\Theta_m}<E_1$ then $\Theta_m$ is linearly stable with respect to $\Theta_n$.
\end{conjecture}


When $\frac{n^2}{m^2}\in I_S$, Theorem \ref{stability12} leaves open the question whether one has stability {\em for any} energy $E_{\Theta_m}>0$,
that is, if a statement similar to Theorem \ref{stability11} holds ($E_2=E_1$?).
According to our numerical computations this does not seem to be the case.
Hence it is reasonable to formulate the following

\begin{conjecture} \label{conj:13}
Assume that $0\le P<m^2<n^2$, $\frac{n^2}{m^2}\in I_S$ and let $E_1$ and $E_2$ be as in Theorem \ref{stability12}. Then there exist $0<E_1\le E_3<E_4\le E_2$ such that:\par
$\bullet$ if $0<E_{\Theta_m}<E_1$ or $E_{\Theta_m}>E_2$ then $\Theta_m$ is linearly stable with respect to $\Theta_n$;\par
$\bullet$ if $E_3<E_{\Theta_m}<E_4$ then $\Theta_m$ is linearly unstable with respect to $\Theta_n$.
\end{conjecture}

\begin{center}
\begin{figure}[ht]
\includegraphics[width=0.4\textwidth]{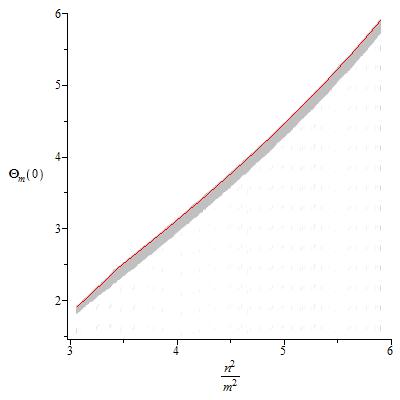}\quad\includegraphics[width=0.4\textwidth]{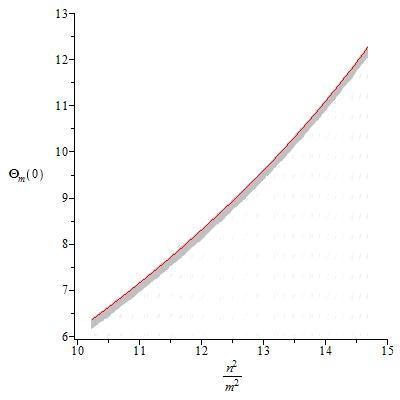}
\caption{An instability region when $\frac{n^2}{m^2} \in (3,6)$ (left) and when $\frac{n^2}{m^2} \in (10,15)$ (right).}
\label{plot1}
\end{figure}
\end{center}

It appears possible that there are many alternating intervals of energy yielding stability or instability, see \cite{bgz}.
In order to clarify this and other questions, we performed our numerical analysis in the special case $P=0$ and $m^2<n^2$. We found a clear evidence for the validity of Conjecture \ref{conj:13}.
More precisely, whenever $\frac{n^2}{m^2} \in I_S$ we found $E_3<E_4$ as in Conjecture \ref{conj:13}.
The natural question then becomes: \emph{is there just one instability region or, equivalently, $E_1=E_3$ and $E_4=E_2$?}
Our numerical experiments prove that, at least generally, this  is not true. For instance setting $m=3$ and $n=7$, namely $\frac{49}{9}\approx 5.44$,
we found that there are at least two instability regions. The instability regions in this case are very narrow and,
for this reason, we only draw one of them in Figure \ref{plot1}. It is known that the geometry of the instability regions may be fairly complicated for general Hill equations, see \cite{broer,broer2}.
\par
We focus now on the case $\frac{n^2}{m^2} \in (3,6)$ and $\frac{n^2}{m^2} \in (10,15)$, which are the
first stability intervals for high energies.
The dependence rule of $E_3$ and $E_4$ with respect to $\frac{n^2}{m^2}$ appears hard to figure out while it is instead more convenient to
consider the dependence of $\Theta_m(0)$ with respect to $\frac{n^2}{m^2}$.
In Figure \ref{plot1} we represent some values of $\Theta_m(0)$ which generate instability for the intervals
$(3,6)$ and $(10,15)$ of $\frac{n^2}{m^2}$. If $\Theta_m(0)$ belongs to the shaded region and
$\dot{\Theta}_m(0)=0$, then the mode $\Theta_m$ is linearly unstable with respect to $\Theta_n$; as remarked above, this is probably not the only instability region.

\begin{center}
\begin{figure}[ht]
\includegraphics[width=0.4\textwidth]{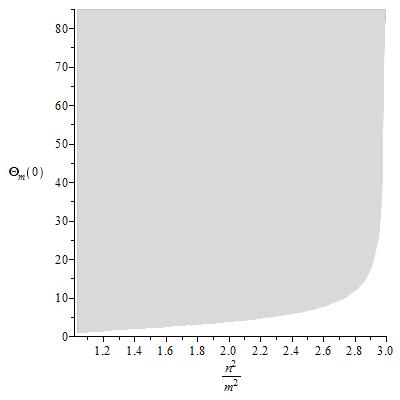}
\quad \includegraphics[width=0.4\textwidth]{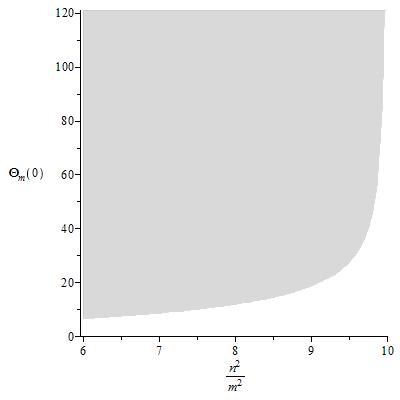}
\caption{Instability regions for high energy when $\frac{n^2}{m^2} \in (1,3)$ (left)
and when $\frac{n^2}{m^2} \in (6,10)$ (right).}
\label{plot:inst}
\end{figure}
\end{center}

Let us now consider the case of instability for high energy, that is, $\frac{n^2}{m^2}\in I_U$.
Analyzing the instability interval $(1,3)$, we found, for several couples
$(n,m)$, $E_3, E_4, E_5,E_6$ such that if
$E_1\leq E_3<E_{\Theta_m}<E_4<E_2$ then $\Theta_m$ is unstable and
if $ E_5<E_{\Theta_m}<E_6\leq E_2$ then $\Theta_m$ is linearly stable.
That is, also in the unstable case it is not true that $E_1=E_2$.
In Figure \ref{plot:inst}, we plot the high energy instability regions for the intervals
$(1,3)$ and $(6,10)$.
 As in Figure \ref{plot1}, we consider the dependence of
 $\Theta_m(0)$ with respect to $\frac{n^2}{m^2}$, if $\Theta_m(0)$ belongs to the shaded region
 and $\dot{\Theta}_m(0)=0$,
 then the mode $\Theta_m$ is linearly unstable with respect to $\Theta_n$.

Theorems \ref{stability12} and \ref{t:limit-case} leave open the
question of stability for large energies when $\frac{n^2}{m^2}\not\in I_U\cup I_S$.
There are infinitely many such couples, for instance,
$(m,n)\in\{(h,6h);(h,35h);(h,204h);(h,1189h);...\, h\in\mathbb{N}\}$.
Even if this is an ``infrequent'' case, it is interesting to notice that numerical results suggest that a different stability behavior occurs
at the endpoints of the intervals $((k+1)(2k+1),(k+1)(2k+3))$. More precisely, we formulate

\begin{conjecture} \label{conj:14}
Assume that $0\le P<m^2<n^2$ and let $k\in \N$. Then there exist $E_k>0$ such that:\par
$\bullet$ if $\frac {n^2}{m^2}= (k+1)(2k+1)$ and $E_{\Theta_m}>E_k$ then $\Theta_m$ is linearly unstable with respect to $\Theta_n$;\par
$\bullet$ if $\frac{n^2}{m^2} = (k+1)(2k+3)$ and $E_{\Theta_m}>E_k$ then $\Theta_m$ is linearly stable with respect to $\Theta_n$.
\end{conjecture}

 The curves that bound the shaded regions in Figure \ref{plot:inst}
 are increasing with respect to $\frac{n^2}{m^2}$,
 Conjecture \ref{conj:14} is related to the fact that they diverge to $+ \infty$ as
 $\frac{n^2}{m^2}$ tends, respectively, to $3$ and $10$.

\subsection{The saddle point case}

\begin{figure}[ht]
\begin{center}
 {\includegraphics[height=44mm, width=76mm]{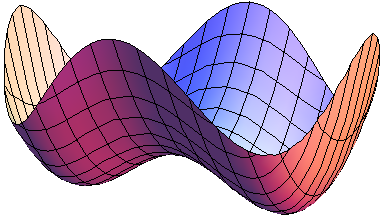}}
\caption{The potential energy functional ${\mathcal U }$ when
$\min\{m^2,n^2\}<P\le \max\{m^2,n^2\}$.}\label{Pinmezzo}
\end{center}
\end{figure}

If $\min\{m^2,n^2\}<P\le \max\{m^2,n^2\}$, then the
functional ${\mathcal U }$ in \eq{potential} is convex in one direction and
has a double well in the orthogonal direction; its qualitative
graph is plotted in Figure \ref{Pinmezzo}. The topology and
geometry of the sublevels of ${\mathcal U }$ depend on the level
considered; this plays an important role in the stability
analysis.

In this situation we have

\begin{theorem}\label{stability2}
If $n^2<P\le m^2$, then there exist $0<E_1\le E_2<\infty$ such that:
\begin{itemize}
\item[(i)] if $0<E_{\Theta_m}<E_1$ then
$\Theta_m$ is linearly unstable with respect to $\Theta_n$;
\item[(ii)] if $E_{\Theta_m}>E_2$ then $\Theta_m$ is linearly stable with respect to $\Theta_n$.
\end{itemize}

\end{theorem}

Again, it would be interesting to understand whether $E_1=E_2$. Note that when $m^2>P$ we necessarily have $E_{\Theta_m}>0$,
see \eq{energies3}. This is no longer true when $m^2<P$ and a different statement holds.

\begin{theorem}\label{stability22}
Assume that $m^2<P\le n^2$, let $I_U$ and $I_S$ as in \eqref{eq:condinst}.
\begin{itemize}
\item[(i)] If either
\neweq{PPP}
L:=\frac{n}{m}\, \sqrt{\frac{2(n^2-m^2)}{P-m^2}}\not\in\N
\endeq
or
\begin{equation}
\label{PPP-2} L \in \N \quad \text{and} \quad 3m^4L^4-(3m^4+4n^2m^2)L^2+4n^2m^2-4n^4\neq 0 \, ,
\end{equation}
there exists $-\frac{(P-m^2)^2}{4}<E_1\le 0$ such that if $-\frac{(P-m^2)^2}{4}<E_{\Theta_m}<E_1$ then $\Theta_m$ is
linearly stable with respect to $\Theta_n$.
\item[(ii)] There exists $E_2\ge0$ such that if $\frac{n^2}{m^2}\in I_U$ and $E_{\Theta_m}>E_2$ then $\Theta_m$
is linearly unstable with respect to $\Theta_n$.
\item[(iii)] There exists $E_2\ge0$ such that if $\frac{n^2}{m^2}\in I_S$ and $E_{\Theta_m}>E_2$ then $\Theta_m$ is
linearly stable with respect to $\Theta_n$.

\end{itemize}
\end{theorem}

Assumptions \eq{PPP} and \eqref{PPP-2} deserve several comments.
Condition \eq{PPP} is needed in order to rule out the resonance
cases in the stability regions of the Hill equation \eq{hill2};
the proof of the stability result under condition \eqref{PPP} is
based on a criterion by Zhukovskii, see Proposition \ref{lyapzhu}
(i). On the other hand, under condition \eqref{PPP-2}, the above
mentioned criterion is no more applicable and the proof is then
based on a refined asymptotic expansion of components of the
monodromy matrix associated to the Hill equation \eqref{hill2} as
the energy $E_{\Theta_m}$ approaches $-\frac{(P-m^2)^2}4$. In
order to avoid vanishing of the higher order term in our
asymptotic expansion, we need the algebraic condition contained in
\eqref{PPP-2}. Otherwise, a higher order asymptotic expansion
should be needed to give an answer to the question of linear
stability of $\Theta_m$. However, we numerically checked that no integer root
$L$ of the fourth order polynomial appearing in \eqref{PPP-2}
exists, at least for $n=2,\dots,5000$ and $m=1,\dots,n-1$.


\subsection{The local maximum case}

\begin{figure}[ht]
\begin{center}
 {\includegraphics[height=44mm, width=76mm]{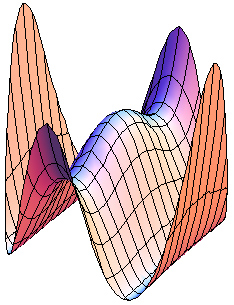}}
\caption{The potential energy functional ${\mathcal U }$ when $\max\{m^2,n^2\}<P$.}\label{Pgrande}
\end{center}
\end{figure}

If $\max\{m^2,n^2\}<P$, then the functional ${\mathcal U }$ in \eq{potential} admits a local maximum at the origin although it remains globally coercive;
its qualitative graph is plotted in Figure \ref{Pgrande}. Also in this case the stability strongly depends on the topology and
geometry of the sublevels of ${\mathcal U }$.

The next result is quite similar to Theorem \ref{stability2}; since also its proof is similar we prove them both in Section \ref{23}.
\begin{theorem}\label{stability3}
Assume that $n^2<m^2<P$. Then there exists $0<E_2<\infty$ such that:
\begin{itemize}
\item[(i)] if $-\frac{(P-m^2)^2}{4}<E_{\Theta_m}\le-\frac{(P-m^2)^2}{4}+\frac{(m^2-n^2)^2}{4}$ and $E_{\Theta_m}\neq0$ then $\Theta_m$ is
linearly unstable with respect to $\Theta_n$; in particular, if $2m^2-n^2-P\ge0$ then the linear instability occurs whenever $E_{\Theta_m}<0$;
\item[(ii)] if $E_{\Theta_m}>E_2$ then $\Theta_m$ is linearly stable with respect to $\Theta_n$.
\end{itemize}
\end{theorem}

Finally, let us also examine the last possible combination of
$m,n,P$. We observe that the statement of the next theorem is
completely similar to the one of Theorem \ref{stability22}. The
only difference consists in the position of $P$ with respect to
$m^2<n^2$.

\begin{theorem}\label{stability4}
Assume that $m^2<n^2<P$, let $I_U$ and $I_S$ as in
\eqref{eq:condinst}.
\begin{itemize}
\item[(i)] If either \eqref{PPP} or \eqref{PPP-2} hold true,
there exists $-\frac{(P-m^2)^2}{4}<E_1\le 0$ such that $\Theta_m$
is linearly stable with respect to $\Theta_n$ provided that
$-\frac{(P-m^2)^2}{4}<E_{\Theta_m}<E_1$.
 \item[(ii)] There exists $E_2\ge0$ such that
if $\frac{n^2}{m^2}\in I_U$ and $E_{\Theta_m}>E_2$ then $\Theta_m$
is linearly unstable with respect to $\Theta_n$.

\item[(iii)] There exists $E_2\ge0$ such that if
$\frac{n^2}{m^2}\in I_S$ and $E_{\Theta_m}>E_2$ then $\Theta_m$ is
linearly stable with respect to $\Theta_n$.

\end{itemize}
\end{theorem}

In Table \ref{riassunto} we summarize all the stability results obtained in the previous statements. It appears that they depend on the order of $P$, $m^2$, $n^2$ and, for large energies, the ratio $\frac{n^2}{m^2}$ is the relevant parameter.

\begin{table}[ht]
\begin{center}
\begin{tabular}{|c|c|c|c|}
\hline
$m,n,P\downarrow$ Energy$\to$ & Low & High & Theorem \\
\hline
$P\le n^2<m^2$ & S & S & \ref{stability0}-\ref{stability11}  \\
$n^2<P\le m^2$ & I & S & \ref{stability2} \\
$n^2<m^2<P$  & I & S & \ref{stability3}  \\
$P<m^2<n^2$  & S & I/S & \ref{stability12} \\
$P=m^2<n^2$ & ? & I/S & \ref{t:limit-case}  \\
$m^2<P\le n^2$  & S & I/S & \ref{stability22}  \\
$m^2<n^2<P$  & S & I/S & \ref{stability4}  \\
\hline
\end{tabular}
\caption{Linear stability (S) and instability (I) of $\Theta_m$ with respect to $\Theta_n$; low means near to the least available energy.}\label{riassunto}
\end{center}
\end{table}

\section{Proof of Theorems \ref{periodicT1} and \ref{periodicT2}}\label{pT1}

The existence, uniqueness and periodicity of the solution $\Theta_k$ of \eq{ODE}-\eq{alphabeta} is a known fact from the theory of ODE's.
By varying the initial data $\alpha$ and $\beta$ we vary $E(\alpha,\beta)$ and obtain infinitely many periodic-in-time solutions of \eqref{truebeam}
in the form \eqref{form}. The solutions in closed form may be expressed in terms of elliptic functions, see \cite{burg}.\par\smallskip
If $k^2\ge P$, for a given $E>0$ we may rewrite \eq{energyy} as
\neweq{ancoraenergia}
\frac{2}{k^4}\dot{\Theta}_k^2=(\Lambda_1-\Theta_k^2)(\Lambda_2+\Theta_k^2)
\endeq
where we omitted the argument $(E)$ of both the $\Lambda_i$'s: note that $\Lambda_2\ge\Lambda_1>0$ (if $k^2=P$ we have
$\Lambda_2=\Lambda_1=2\sqrt{E}/k^2$). Hence,
$$
\|\Theta_k\|_\infty=\sqrt{{\Lambda_1}}\, ,\qquad-\sqrt{{\Lambda_1}}\le\Theta_k(t)\le\sqrt{{\Lambda_1}}\quad\forall t
$$
and $\Theta_k$ oscillates in this range. If $\Theta_k(t)$ solves \eq{ODE}-\eq{alphabeta} for $\beta=0$, then also $\Theta_k(-t)$ solves the same problem:
this shows that the period $T(E)$ of a solution $\Theta_k$ of \eq{ODE} is the double of the length of an interval of monotonicity for $\Theta_k$.
Since the problem is autonomous, we may assume that $\Theta_k(0)=-\sqrt{\Lambda_1}$ and $\dot{\Theta}_k(0)=0$; then we have that
$\Theta_k(T(E)/2)=\sqrt{\Lambda_1}$ and $\dot{\Theta}_k(T(E)/2)=0$. By rewriting \eq{ancoraenergia} as
$$
\frac{\sqrt2}{k^2}\dot{\Theta}_k=\sqrt{(\Lambda_1-\Theta_k^2)(\Lambda_2+\Theta_k^2)}\qquad\forall t\in\left(0,\frac{T(E)}{2}\right)\, ,
$$
by separating variables, and upon integration over the time interval $(0,T(E)/2)$ we obtain
$$\frac{T(E)}{2}=\frac{\sqrt2}{k^2}\int_{-\sqrt{\Lambda_1}}^{\sqrt{\Lambda_1}}\frac{d\theta}{\sqrt{(\Lambda_2+\theta^2)(\Lambda_1-\theta^2)}}\,.$$
Then, using the fact that the integrand is even with respect to $\theta$ and through a change of variable, we obtain \eq{TE}.\par
Both the maps $E \mapsto \Lambda_i(E)$ are continuous and increasing for $E\in[0,\infty)$ and $\Lambda_1(0)=0$, $\Lambda_2(0)=2(1-P/k^2)$.
Whence, $E\mapsto T(E)$ is strictly decreasing and \eq{limitperiod} holds; if $k^2>P$ this limit could have also been obtained by noticing that,
as $E\downarrow0$, the equation \eq{ODE} ``tends'' to the linear equation $\ddot{\Theta}_k(t)+k^2(k^2-P)\Theta_k(t)=0$. This completes the proof of Theorem
\ref{periodicT1} (case $k^2\ge P$).\par\smallskip
If $k^2<P$ and $E(\alpha,\beta)>0$, then $\Lambda_1(E)>\Lambda_2(E)>0$. The same arguments as above yield that the period of $\Theta_k$ is given by \eq{TE}.
Since both $\Lambda_{1,2}(E)\to\infty$ as $E\to\infty$, we infer the last limit in \eq{limitperiod2}. Moreover, here we have $\Lambda_2(0)=0$
and $\Lambda_1(0)=2(\frac{P}{k^2}-1)>0$, which proves the limit value \eq{limitperiod2} as $E\downarrow0$.\par
If $k^2<P$ and $-\frac{(P-k^2)^2}{4}<E(\alpha,\beta)<0$, we set
$$\Phi_1(E):=\frac{P-k^2+\sqrt{(P-k^2)^2+4E}}{k^2}>0\, ,\quad\Phi_2(E):=\frac{P-k^2-\sqrt{(P-k^2)^2+4E}}{k^2}>0$$
and we notice that $\Phi_1>\Phi_2$. Then, instead of \eq{ancoraenergia}, we obtain
\neweq{duedue}
\frac{2}{k^4}\dot{\Theta}_k^2=(\Phi_1-\Theta_k^2)(\Theta_k^2-\Phi_2)\, .
\endeq
This readily shows that $\Phi_2\le\Theta_k(t)^2\le\Phi_1$ for all $t$ and $\Theta_k^2$ oscillates in this range. Let us assume that $\Theta_k>0$,
since the case $\Theta_k<0$ is completely similar: $0<\sqrt{\Phi_2}\le\Theta_k(t)\le\sqrt{\Phi_1}$ for all $t$.\par
Again, the period $T(E)$ of a solution $\Theta_k$ of \eq{ODE} is the double of the length of an interval of monotonicity for $\Theta_k$.
We take $\Theta_k(0)=\sqrt{\Phi_2}$ and $\dot{\Theta}_k(0)=0$; then we have that $\Theta_k(T(E)/2)=\sqrt{\Phi_1}$ and $\dot{\Theta}_k(T(E)/2)=0$.
By rewriting \eq{duedue} as
$$
\frac{\sqrt2}{k^2}\dot{\Theta}_k=\sqrt{(\Phi_1-\Theta_k^2)(\Theta_k^2-\Phi_2)}\qquad\forall t\in\left(0,\frac{T(E)}{2}\right)\, ,
$$
by separating variables, and upon integration over the time interval $(0,T(E)/2)$ we obtain
$$\frac{T(E)}{2}=\frac{\sqrt2}{k^2}\int_{\sqrt{\Phi_2}}^{\sqrt{\Phi_1}}\frac{d\theta}{\sqrt{(\Phi_1-\theta^2)(\theta^2-\Phi_2)}}\,.$$
Then, after a change of variable and replacing $\Phi_1$, $\Phi_2$, and $\delta$, we get \eq{TE2}.\par
If $E\uparrow0$ then $\delta\to0$ and $T(E)\to+\infty$. If $E\downarrow-\frac{(P-k^2)^2}{4}$ then $\delta\to1$; moreover, by recalling that
$$\int_\delta^1\frac{ds}{\sqrt{(1-s)(s-\delta)}}=\pi\qquad\forall\delta\in(0,1)\, ,$$
we obtain the estimates
$$
\frac{\pi}{\sqrt{2(1+\delta)}}\le\int_\delta^1\frac{d\theta}{\sqrt{(1-\theta^2)(\theta^2-\delta^2)}}\le\frac{\pi}{\sqrt{2\delta(1+\delta)}}
$$
which prove the first limit in \eq{limitperiod2} by letting $\delta\to1$.\par
Finally, the monotonicity of the period $T=T(E)$ follows from \eq{TE} when $E>0$ and from a result by Chicone \cite{Ci87} (see also
\cite[Theorem 1]{yaga}) when $E<0$. This completes the proof of Theorem \ref{periodicT2} (case $k^2<P$).

\section{Bounds for the amplitudes and periods}\label{tech}

Consider the equation \eq{ODE} with $k=m$ and take initial conditions with no kinetic energy:
\neweq{upto}
\ddot{\Theta}_m(t)+m^2(m^2-P)\Theta_m(t)+m^4\Theta_m(t)^3=0\, ,\qquad
\Theta_m(0)>0\, ,\quad\dot{\Theta}_m(0)=0\, .
\endeq

\begin{lemma}\label{stimaw}
Let $P\in\R$, $m\in\N$, let $\Theta_m$ be the solution of
\eqref{upto} and let $E>0$ denote its energy.
\begin{itemize}
\item[(i)] If $P<m^2$
$$
\|\Theta_m\|_\infty^2=\frac{2}{m^2(m^2-P)}\, E+o(E)\, ,\qquad\mbox{as }E\to0\, .
$$

\item[(ii)] For any $P\in\R$ and $m\in\N$ we have
$$
\|\Theta_m\|_\infty^2=\frac{2}{m^2}\, \sqrt{E}+o\big(\sqrt{E}\big)\qquad\mbox{as }E\to\infty\, .
$$
\end{itemize}
\end{lemma}
{\em Proof.} Since $\dot{\Theta}_m(0)=0$ we have
\neweq{Econ}
E=\frac{m^2(m^2-P)}{2}\, \Theta_m(0)^2+\frac{m^4}{4}\, \Theta_m(0)^4\qquad\mbox{and}\qquad
\|\Theta_m\|_\infty^2=\frac{\sqrt{(m^2-P)^2+4E}-m^2+P}{m^2}\, .
\endeq
From \eq{Econ} we infer the two asymptotic estimates.\endproof

Let us introduce a constant which will be of great importance in the sequel:
\neweq{sigma}
\sigma:=\int_0^1\frac{d\theta}{\sqrt{1-\theta^4}}=\frac{1}{\sqrt2}K\big(\tfrac{1}{\sqrt2}\big)\approx1.311\, ,
\endeq
where $K(\cdot)$ denotes the complete elliptic integral of the first kind, that is,
\begin{equation*}
K(x)=\int_0^1 \frac{1}{\sqrt{(1-t^2)(1-x^2t^2)}}\, dt=\int_0^{\frac \pi 2} \frac 1{\sqrt{1-x^2 \sin^2\alpha}}\, d\alpha \qquad \text{for any } x\in (0,1) \, .
\end{equation*}
The representation of $\sigma$ in terms of $K$, follows from the change of variables $\theta=\sqrt{1-t^2}$.\par
We now prove some asymptotic estimates when the energy $E$ tends to both $0$ and $\infty$.

\begin{lemma}
Let $P\in\R$, $m\in\N$, let $\Theta_m$ be the solution of
\eqref{upto}, let $E>0$ denote its energy and let $T(E)$ denote
its period.
\begin{itemize}
\item[(i)] If $P<m^2$
\neweq{stimavera}
\left(\frac{2\pi}{T(E)}\right)^2=m^2(m^2-P)+\frac32
\frac{m^2}{m^2-P}\, E+o(E)\, ,\qquad\mbox{as }E\to0\, .
\endeq

\item[(ii)] For any $P\in\R$ and $m\in\N$ we have
\neweq{stimavera2}
T(E)=\frac{4\sigma}{m\sqrt[4]{E}}+o(1/\sqrt[4]{E})\qquad\mbox{as }E\to\infty\, ,
\endeq
where $\sigma$ is as in \eqref{sigma}.
\end{itemize}
\end{lemma}

{\em Proof.} From \eq{Lambdas} we infer that
$$
\Lambda_1(E)=\frac{2}{m^2(m^2-P)}\, E+o(E)\, ,\quad\Lambda_2(E)
=\frac{2(m^2-P)}{m^2}+\frac{2}{m^2(m^2-P)}\, E+o(E)\,
,\qquad\mbox{as }E\to0\, .
$$
By plugging these estimates into \eq{TE} and with some tedious computations we obtain
$$
T(E)=\frac{2\pi}{m\sqrt{m^2-P}}-\frac{3\pi}{2m(m^2-P)^{5/2}}\, E+o(E)\, ,\qquad\mbox{as }E\to0\, .
$$
From this estimate we then obtain \eq{stimavera}. From \eq{Lambdas} and \eq{TE} we obtain \eq{stimavera2}.\endproof

We put
\neweq{IE}
I(E):=\int_0^{T(E)/2}\Big(n^2(n^2-P)+m^2n^2\Theta_m(t)^2\Big)^2\, dt\, .
\endeq
By combining the above asymptotic estimates with energy arguments, we can prove the following result.

\begin{lemma}\label{nuovastima}
Let $\sigma$ be as in \eqref{sigma} and let $I(E)$ be as in \eqref{IE}. For all $n,m\in\N$ and $P\in \R$ we have
\neweq{this}
\lim_{E\to\infty}\left(\frac{T(E)}{2}\right)^3\, I(E)=\frac{64\, n^4}{3\, m^4}\, \sigma^4\, .
\endeq
\end{lemma}
{\em Proof.} Let us first translate the function $\Theta_m$, solution of \eq{upto}, in such a way that $\Theta_m(0)=0$ and $\dot{\Theta}_m(0)>0$;
this also implies that $\Theta_m(T(E)/2)=0$. Therefore, if we multiply the equation in \eq{upto} by $\Theta_m$ and we integrate by parts over $(0,T(E)/2)$
we obtain
\neweq{comb1}
m^2(m^2-P)\int_0^{T(E)/2}\Theta_m(t)^2\, dt+m^4\int_0^{T(E)/2}\Theta_m(t)^4\, dt=\int_0^{T(E)/2}\dot{\Theta}_m(t)^2\, dt\, .
\endeq
On the other hand, by integrating over $(0,T(E)/2)$ the conservation of the energy law \eq{energyy}, we obtain
$$
\int_0^{T(E)/2}\dot{\Theta}_m(t)^2\, dt+m^2(m^2-P)\int_0^{T(E)/2}\Theta_m(t)^2\, dt+\frac{m^4}{2}\int_0^{T(E)/2}\Theta_m(t)^4\, dt=E\cdot T(E)\, .
$$
By combining this with \eq{comb1} we infer that
\neweq{identity}
2m^2(m^2-P)\int_0^{T(E)/2}\Theta_m(t)^2\, dt+\frac{3m^4}{2}\int_0^{T(E)/2}\Theta_m(t)^4\, dt=E\cdot T(E)\, .
\endeq

By \eq{stimavera2} and Lemma \ref{stimaw} we obtain
$$
E\cdot T(E)\sim\frac{4\sigma}{m}\, E^{3/4}\ ,\qquad\int_0^{T(E)/2}\Theta_m(t)^2\, dt\le\frac{\|\Theta_m\|_\infty^2\, T(E)}{2}\sim
\frac{4\sigma}{m^3}\, E^{1/4}=o\big(E^{3/4}\big)\qquad\mbox{as }E\to\infty\, .
$$
By taking this into account, \eq{identity} yields
$$
\int_0^{T(E)/2}\Theta_m(t)^4\, dt\sim\frac{8\sigma}{3m^5}\, E^{3/4}\qquad\mbox{as }E\to\infty\, .
$$
Therefore,
$$I(E)\sim m^4n^4\int_0^{T(E)/2}\Theta_m(t)^4\, dt\sim
\frac{8n^4\sigma}{3m}\, E^{3/4}\qquad\mbox{as }E\to\infty\, .$$
After multiplication by $(T(E)/2)^3$ and using again \eq{stimavera2} we obtain \eq{this}.\endproof

We now introduce three functions of $E$ that will allow to simplify some notations and estimates in the sequel. We define
\neweq{Zm}
X_m(E):=4E+(m^2-P)^2\, ,\quad Y_m(E):=\frac{X_m(E)}{(m^2-P)^2}\, ,\quad Z_m(E):=\frac12- \frac{1}{2\sqrt{Y_m(E)}}\, ;
\endeq
note that these three functions are all strictly increasing with respect to $E$ and that, since $0<E<\infty$, we have the bounds
$$
(m^2-P)^2<X_m(E)<\infty\, ,\quad1<Y_m(E)<\infty\, ,\quad0<Z_m(E)<\frac12\, .
$$

With the change of variables $\theta=\cos\alpha$ and recalling both \eq{Lambdas} and \eq{TE}, we obtain
\neweq{newTE}
T(E)=\frac{4}{m\sqrt[4]{X_m(E)}}\int_0^{\pi/2}\frac{d\alpha}{\sqrt{1-(\frac12- \frac{m^2-P}{2\sqrt{X_m(E)}})\sin^2\alpha}}=
\frac{4}{m\sqrt[4]{X_m(E)}}\, K\left(\sqrt{\frac12- \frac{m^2-P}{2\sqrt{X_m(E)}}}\right)\, ,
\endeq
where $K(\cdot)$ is the complete elliptic integral of the first kind and $X_m$ is as in \eq{Zm}.\par
We conclude this section with a bound for the $L^2$-norm of the solution of \eq{upto}.
From \eq{identity} and the H\"older inequality we infer that
$$
\frac{3m^4}{T(E)}\left(\int_0^{T(E)/2}\Theta_m(t)^2\, dt\right)^2+
2m^2(m^2-P)\int_0^{T(E)/2}\Theta_m(t)^2\, dt-E\cdot T(E)\le0\, .
$$
By solving this second order algebraic inequality we obtain
\neweq{stimaL2}
\int_0^{T(E)/2}\Theta_m(t)^2\, dt\le\frac{T(E)}{3m^2}\Big(\sqrt{3E+(m^2-P)^2}-(m^2-P)\Big)<\frac{T(E)}{3}\frac{\sqrt{X_m(E)}-m^2+P}{m^2}\, .
\endeq

\section{The Cazenave-Weissler result for large energies} \label{s:cw}

The purpose of the present section is to prove

\begin{theorem}\label{stability121}
Let $m$ and $n$ be two positive integers and let $I_U$ and $I_S$ be as in \eqref{eq:condinst}. If $\frac{n^2}{m^2}\in I_U$ (resp.\ $\frac{n^2}{m^2}\in I_S$) there exists $\overline{E}>0$ such that if $E_{\Theta_m}>\overline{E}$ then $\Theta_m$ is linearly unstable (resp.\ stable) with respect to $\Theta_n$.
\end{theorem}

By Definition \ref{defstabb}, we have to analyze the stability of the Hill equation \eqref{hill2}. To simplify the notation we rewrite
system \eqref{cw} as \eq{cwx}; then, the substitution $w(t)\mapsto w\pt{\frac{t}{m}}$ and $z(t)\mapsto z\pt{\frac{t}{m}}$ leads to
\begin{equation*}
\left\{ \begin{array}{l}
\ddot{w}(t)+ (m^2-P) w(t)+( w(t)^2+ z(t)^2)w(t)=0\\
\ddot{z}(t)+ \frac{n^2}{m^2}(n^2-P) z(t)+\frac{n^2}{m^2}( w(t)^2+ z(t)^2)z(t)=0.
\end{array}\right.
\end{equation*}

Setting $\nu=(m^2-P)$, $\nu'=(n^2-P)$ and $\gamma=\frac{n^2}{m^2}$, we obtain
\begin{equation}\label{eq:sitpoinc}
\left\{ \begin{array}{l}
\ddot{w}(t)+ \pt{\nu+ w(t)^2+ z(t)^2} w(t)=0\\
\ddot{z}(t)+ \gamma\pt{\nu' + w(t)^2+ z(t)^2}z(t)=0.
\end{array}\right.
\end{equation}
System \eq{eq:sitpoinc} is Hamiltonian and has conserved energy given by
\begin{align*}
\E\Big(w(t),\dot{w}(t),z(t),\dot{z}(t)\Big)=\frac{\dot{w}(t)^2}{2}+\frac{\dot{z}(t)^2}{2\gamma}+\nu\, \frac{w(t)^2}{2}+\nu'\, \frac{z(t)^2}{2}
+\frac{\big(w(t)^2+z(t)^2\big)^2}{4}\equiv E_0
\end{align*}
for some $E_0$ depending on the initial data.
Then, we may rephrase Theorem \ref{stability121} as follows.

\begin{proposition}\label{rephrase}
 Let $I_U$ and $I_S$ be as in \eqref{eq:condinst}.
 If $\gamma\in I_U$ (resp.\ $\gamma\in I_S$), then there exists $\overline{E}>0$
 such that if $E_0>\overline{E}$ then $w$ is linearly unstable (resp.\ stable) with respect to $z$.
\end{proposition}

The proof of Proposition \ref{rephrase} is essentially due to Cazenave-Weissler \cite{cazw2}, see also \cite{cazw}.
The main idea is to determine the eigenvalues of the Jacobian at the origin of the Poincar\'e map:
if the eigenvalues $\lambda_1=\lambda$, $\lambda_2=\lambda^{-1}$, are real with $|\lambda|\in (0,1)$,
then we have linear instability, if $\lambda_1= e^{i\omega}$, $\lambda_2=e^{-i \omega}$, with
$\omega \in (0,\pi)$ we have linear stability, see \cite[Section 2.4.4]{chicone}. Actually, when $\gamma\in I_U$ and $E_0>\overline{E}$, we also have that $w$ is orbitally unstable with respect to $z$, see the last two lines of the proof.

We briefly sketch the proof by
emphasizing the differences with \cite{cazw2}. In particular, here $\nu\neq\nu'$ and $\nu$ and $\nu'$ may be negative.\par
We fix some $E_0>0$ and we define the following open neighborhood of $(0,0)\in\R^2$:
$$\mathcal{ V}_1=\left\{(a,b)\in\R^2:\, \frac{b^2}{2\gamma}+\nu'\frac{a^2}{2}+\frac{a^4}{4}<E_0\right\}\, ;$$
whence, $\E(0,0,a,b)<E_0$ for all $(a,b)\in\mathcal{ V}_1$. For a given couple $(a,b)\in\mathcal{ V}_1$, we consider the solution $(w,z)$ of \eq{eq:sitpoinc}
with initial conditions
\begin{align}\label{eq:incond}
w(0)=0, \quad \dot{w}(0)=w_1>0, \quad z(0)=a, \quad \dot{z}(0)=b,
\end{align}
where $w_1$ is chosen in such a way that $\E(0,w_1,a,b)=E_0$, that is,
\neweq{w1}
w_1=\sqrt{2E_0-\frac{b^2}{\gamma}-\nu'a^2-\frac{a^4}{2}}.
\endeq
Let us prove the following statement.

\begin{lemma}\label{firstzero}
There exists a neighborhood $\mathcal{ V}_2$ of $(0,0)$ such that if $(a,b)\in\mathcal{ V}_1\cap\mathcal{ V}_2$ and $w_1$ is as in \eqref{w1}, then
there exists a first point $\tau=\tau(a,b)>0$ such that $w(\tau)=0$, where $(w,z)$ is the solution of \eqref{eq:sitpoinc}-\eqref{eq:incond}.
\end{lemma}
{\em Proof.} If $\nu>0$ then there exist infinitely many points $\tau_n$ such that $w(\tau_n)=0$ and $|\tau_{n+1}-\tau_{n}|\leq \frac{\pi}{\sqrt{\nu}}$;
to see this, it suffices to multiply the first equation in \eqref{eq:sitpoinc} by $\sin\pt{ \sqrt{\nu} t}$
and integrate twice by parts on the interval $(0,\frac{\pi}{\sqrt{\nu}})$, see \cite{cazw2}.\par
If $\nu\le0$, this simple trick does not work and we use an abstract argument. We first notice that if $W_\nu$ solves the problem
\begin{align}
  \label{eq:sitpoinclin}
  \left\{ \begin{array}{l}
	\ddot{W}_\nu(t)+  \pt{\nu+ W_\nu(t)^2} W_\nu(t)=0\\
	W_{\nu}(0)=0, \quad \dot{W}_{\nu}(0)=\sqrt{2 E_0}\, ,
	\end{array}\right.
\end{align}
then the couple $(w,z)=(W_\nu,0)$ is a periodic solution of
\eqref{eq:sitpoinc} and $w$ is sign-changing, see Section \ref{nonlinfor}. Then we notice that,
by energy conservation, any solution $(w,z)$ of \eq{eq:sitpoinc}
is globally defined. Hence, by continuous dependence,
there exists a 3-dimensional neighborhood $V_2$ of $(0,\sqrt{2E_0},0,0)$, contained in the hyperplane $x_1=0$, and a map
$$
\tau:\{(x_1,x_2,x_3,x_4):(x_3,x_4)\in \mathcal { V}_1, \,
\mathcal{E}(x_1,x_2,x_3,x_4)=E_0\} \cap V_2 \to \R\ ,\qquad
(0,w_1,a,b)\mapsto \tau(0,w_1,a,b)
$$
such that the solution $w$ of \eqref{eq:sitpoinc}, with initial
condition \eqref{eq:incond}, satisfies $w(\tau(0,w_1,a,b))=0$.
Then one can define $\mathcal{ V}_2$ as the projection of $ V_2$
with respect to the map $(x_1,x_2,x_3,x_4)\mapsto (x_3,x_4)$. This
completes the proof also in the case where $\nu\le0$.\endproof

Note that in Lemma \ref{firstzero} it is essential that
$\dot{w}(0)>0$ since otherwise $w$ could be a one-sign solution of
\eq{eq:sitpoinc} (including constants); this may happen whenever
$\nu<0$. In the sequel, we denote
$$\mathcal{ V}=\mathcal{ V}_1\cap\mathcal{ V}_2\, .$$
Lemma \ref{firstzero} defines the Poincar\'e map
\begin{align*}
T:\mathcal{ V} \to \R^2\ ,\quad T(a,b):=-(z(\tau), \dot{z}(\tau))\ ,\quad\mbox{where }\tau=\tau(0,w_1,a,b)\quad\forall(a,b)\in\mathcal{ V}.
\end{align*}
The map $T$ is defined for all $\nu,\nu',\gamma$; by construction we also know that $T$ is $C^1$ and $T(0,0)=(0,0)$.\par

Let $W_\nu$ be as in \eq{eq:sitpoinclin} and consider the Hill equation
\begin{align}
  \label{eq:hill2}
  \ddot{\xi}(t)+ \gamma\pt{ \nu' + W_\nu(t)^2 }\xi(t)=0,
\end{align}
with initial data $\xi(0)=a$ and $\dot{\xi}(0)=b$. Then we define the linear operator $L:\R^2\to\R^2$ by
\begin{equation}\label{L}
L(a,b)=-\big(\xi(\rho),\dot{\xi}(\rho)\big)\qquad\forall(a,b)\in\R^2
\end{equation}
where $\rho=\rho(E_0)$ is the first positive zero of $W_\nu$ and, in turn, the period of the function $W_\nu(t)^2$, see Theorems \ref{periodicT1} and \ref{periodicT2}. The eigenvalues of $L$ coincide with those of the monodromy matrix of the Hill equation \eq{eq:hill2}, see \cite[Chapter
II, Section 2.1]{yakubovich}.  As in \cite[Proposition 2.1]{cazw2}, one can prove that the Jacobian of $T$ at the origin
coincides with $L$, namely $DT(0,0)=L$.\par
Next, we consider the solution $u$ of the problem
\begin{align}
	\label{eq:u}
	\ddot{u}(t)+  u(t)^3 =0\, ,\qquad u(0)=0\, , \ \dot{u}(0)=1\, .
\end{align}
Then $u$ is a periodic function and changes sign infinitely many times. Denote by $\theta>0$ the first positive zero of $u$ and by $\eta$ the solution of the Hill equation
\begin{align*}
	\ddot{\eta}(t)+\gamma u(t)^2\eta(t)=0\, ,\qquad \eta(0)=a\, ,\ \dot{\eta}(0)=b\, .
\end{align*}
Finally, we define the map
\[
	B_\gamma:\R^2\to\R^2\, ,\qquad B_\gamma(a,b)=-\big(\eta(\theta),\dot{\eta}(\theta)\big)\, .
\]
Let $I_U$ and $I_S$ be as in \eq{eq:condinst}. By \cite[Theorem 3.1]{cazw2} we know that:
\begin{itemize}
	\item if $\gamma \in I_U$ then $B_\gamma$ has eigenvalues $\lambda,\lambda^{-1}\in\R$, for some $|\lambda|\in(0,1)$;
	\item if $\gamma \in I_S$ then $B_\gamma$ has eigenvalues $e^{i\omega}, e^{-i \omega}$, for some $\omega\in(0,\pi)$.
\end{itemize}

Then we perturb \eq{eq:u}.
Let $\nu\neq0$ and, for all $\eps\neq 0$ having the same sign as $\nu$, consider the solution $u_\eps$ of
\[
	\ddot{u}_\eps(t)+\eps u_\eps(t)+u_\eps(t)^3 =0\, ,\qquad u_\eps(0)=0\, , \ \dot{u}_\eps(0)=1\, .
\]
Denote by $\theta_\eps>0$ the first positive zero of $u_\eps$ and consider the problem
\begin{align}
	\label{eq:nudiv0}
	\ddot{\eta}_\eps(t)+\gamma\left(\eps \frac{\nu'}{\nu}+u_\eps(t)^2 \right)\eta_\eps(t)=0\, ,\qquad \eta_\eps(0)=a\, ,\ \dot{\eta}_\eps(0)=b\,
\end{align}
whose solution $\eta_\eps$ defines the map
\[
  B_{\gamma,\eps}:\R^2\to\R^2\, ,\qquad B_{\gamma,\eps}(a,b)=-\big(\eta_\eps(\theta_\eps),\dot{\eta}_\eps(\theta_\eps)\big)\, .
\]
Then $B_{\gamma,\eps}\to B_\gamma$ as $\eps \to 0$.
Therefore, by the above statements and a continuity argument, there exists
$\eps_0>0$ such that if $|\eps|<\eps_0$ there holds
\begin{itemize}
\item if $\gamma\in I_U$ then $B_{\gamma,\eps}$ has eigenvalues $\lambda_\eps, \lambda_\eps^{-1}\in \R$, for some
 $|\lambda_\eps| \in (0,1)$;
\item if $\gamma\in I_S$ then $B_{\gamma,\eps}$ has eigenvalues $e^{i\omega_\eps}$ and
 $e^{-i\omega_\eps}$, for some $\omega_\eps\in(0,\pi)$.
\end{itemize}
 If $W_\nu(t)$ solves \eqref{eq:sitpoinclin} with $E_0=\frac{\nu^2}{2 \varepsilon^2}$ and $\xi(t)$ solves \eqref{eq:hill2} with $\xi(0)=a$ and $\dot{\xi}(0)=b \sqrt{\frac{\nu}{\varepsilon}}$, then
   \[
     W_\nu(t)=\sqrt{\tfrac{\nu}{\varepsilon}} \, u_\varepsilon\pt{\sqrt{\tfrac{\nu}{\varepsilon}}\,t},
     \qquad \xi(t)= \eta_\varepsilon\pt{\sqrt{\tfrac{\nu}{\varepsilon}}\,t}.
   \]
   By direct computation, one checks that the eigenvalues of $B_{\gamma,\eps}$ are the same of $L$ (defined in \eq{L})
   with energy $E_0= \frac{\nu^2}{2 \varepsilon^2}$.
   Therefore, when $E_0>\frac{\nu^2}{2\eps_0}$,
   if $\gamma \in I_U$ then  the system \eqref{eq:sitpoinc} is linearly unstable, while if $\gamma \in I_S$ then the system \eqref{eq:sitpoinc} is linearly stable.

  If $\nu=0$, we replace \eqref{eq:nudiv0} with
  \[
	\ddot{\eta}_\eps(t)+\gamma\left(\eps \nu'+u(t)^2 \right)\eta_\eps(t)=0\, ,\qquad \eta_\eps(0)=a\, ,\  \dot{\eta}_\eps(0)=b\,
  \]
 where $\eps>0$ and $u$ is the solution of \eqref{eq:u}. Furthermore, in this case we have
  \[
     W_\nu(t)=\tfrac{1}{\sqrt{\varepsilon}}\, u \pt{\tfrac{t}{\sqrt{\varepsilon}}},
     \qquad \xi(t)= \eta_\varepsilon\pt{\tfrac{t}{\sqrt{\varepsilon}}}.
   \]
   Proceeding as in the case $\nu\neq 0$, we reach the same conclusion on linear stability and linear instability.

 At this point, in order to prove the orbital instability when $\gamma\in I_U$ and $E_0$ is large enough, one may proceed as in the proof of Theorem 1.1 and Theorem 2.2 in \cite{cazw2}.

\section{Proof of Theorem \ref{stability0}}

If $T=T(E)$ denotes the period of $\Theta_m$, then the function
$a(t)$ in \eq{hill2} has period $T(E)/2$. Since $m^2>n^2\ge P$, by
\eq{limitperiod} we know that
$$\lim_{E\to0}\Big(n^2(n^2-P)+n^2m^2\|\Theta_m\|_\infty^2\Big)=n^2(n^2-P)<m^2(m^2-P)=\lim_{E\to0}\frac{4\pi^2}{T(E)^2}\, .$$
Hence, by continuity, there exists $E_1>0$ such that
$$n^2(n^2-P)+n^2m^2\|\Theta_m\|_\infty^2\le\frac{4\pi^2}{T(E)^2}\qquad\forall E\le E_1\, .$$
Then the first criterion in Proposition \ref{lyapzhu} (with
$\ell=0$) ensures that the trivial solution of \eq{hill2} is
stable, provided that $E\le E_1$.\par Since $m>n$, \eq{this}
proves that
$$
\left(\frac{T(E)}{2}\right)^3I(E)<\frac{64}{3}\, \sigma^4
$$
for sufficiently large $E$. Then the second criterion in
Proposition \ref{lyapzhu} ensures that the trivial solution of
\eq{hill2} is stable, provided that $E$ is sufficiently large, say
for $E>E_2$.\par What we have seen proves the linear stability of
$\Theta_m$ for both $E\le E_1$ and $E>E_2$; Theorem
\ref{stability0} is so proved.

\section{Proof of Theorem \ref{stability11}}

For our convenience we put $\eps=\tfrac{21}{22}$. Since $\eps^2m^2\ge n^2$ and $P\ge0$, we also know that
\neweq{firstcase}
\frac{n^4}{m^4}\le\eps^4\ ,\quad(n^2-P)\le\eps^2(m^2-P)\ ,\quad\quad(n^2-P)^2\le\eps^4(m^2-P)^2\, .
\endeq
Let $\Theta_m$ be the solution of \eq{upto}: denote by $E$ its energy, see \eq{Econ}, and by $T(E)$ its period, see \eq{TE}.\par
We first prove an important implication.

\begin{lemma}\label{bjorn1}
Let $\eps=\tfrac{21}{22}$ and assume that \eqref{firstcase} holds. Let $Z_m$ be as in \eqref{Zm} and let $I(E)$ be as in \eqref{IE}.
Then the following implication holds:
$$K\Big(\sqrt{Z_m(E)}\Big)^4\le\frac{4\, \sigma^4}{\eps^4\, \Big(4(3\eps^4-4\eps^2+\tfrac53)Z_m(E)^2-4(3\eps^4-2\eps^2+\tfrac13)Z_m(E)+3\eps^4\Big)}$$
$$\Longrightarrow\quad\left(\frac{T(E)}{2}\right)^3 I(E)<\frac{64}{3}\, \sigma^4\, .$$
\end{lemma}
{\em Proof.} By computing the squared integrand in \eq{IE} we obtain the bound
\begin{eqnarray*}
I(E) &=& n^4(n^2-P)^2\, \frac{T(E)}{2}+2m^2n^4(n^2-P)\int_0^{T(E)/2}\Theta_m(t)^2\, dt+m^4n^4\int_0^{T(E)/2}\Theta_m(t)^4\, dt\\
\mbox{by \eq{identity}} &=& n^4\, \left((n^2-P)^2+\frac43 E\right)\, \frac{T(E)}{2}\, +
\frac23 \, m^2n^4\Big(3(n^2-P)-2(m^2-P)\Big)\int_0^{T(E)/2}\Theta_m(t)^2\, dt\\
\mbox{by \eq{firstcase} } &\le& \frac{n^4}{3}\Big((3\eps^4\!-\!1)(m^2\!-\!P)^2+X_m(E)\Big)\frac{T(E)}{2}+
\frac23 m^2n^4(3\eps^2\!-\!2)(m^2\!-\!P)\int_0^{T(E)/2}\Theta_m(t)^2dt\\
\mbox{by \eq{stimaL2} } &<& \frac{n^4}{3}\, \frac{T(E)}{2}\, \left(X_m(E)+\tfrac{4}{3}(3\eps^2-2)(m^2-P)\,
\sqrt{X_m(E)}+\big(3\eps^4-4\eps^2+\tfrac53\big)(m^2-P)^2\right)
\end{eqnarray*}
where $X_m$ is as in \eq{Zm}; if we take $Y_m$ and $Z_m$ as in \eq{Zm}, then the latter inequality and \eq{newTE} yield
\begin{eqnarray*}
\left(\frac{T(E)}{2}\right)^3 I(E) &<& \frac{16}{3}\, \frac{n^4}{m^4}\,
\left(1+\frac{4(3\eps^2-2)}{3\, \sqrt{Y_m(E)}}+\frac{3\eps^4-4\eps^2+\tfrac53}{Y_m(E)}\right)\,
K\left(\sqrt{\frac12- \frac{1}{2\sqrt{Y_m(E)}}}\right)^4\\
\mbox{by \eq{firstcase} } &\le& \frac{16}{3}\, \eps^4\, \Big(4(3\eps^4\!-\!4\eps^2\!+\!\tfrac53)Z_m(E)^2-
4(3\eps^4\!-\!2\eps^2\!+\!\tfrac13)Z_m(E)+3\eps^4\Big)\, K\left(\sqrt{Z_m(E)}\right)^4\, .
\end{eqnarray*}
The statement is so proved.\endproof

The second step is another crucial implication.

\begin{lemma}\label{bjorn}
Let $\eps=\tfrac{21}{22}$ and assume that \eqref{firstcase} holds. Let $Z_m$ be as in \eqref{Zm}.
Then the following implication holds:
$$K\Big(\sqrt{Z_m(E)}\Big)^2\le\frac{\pi^2}{4\eps^2\Big(\eps^2+2(1\!-\!\eps^2)Z_m(E)\Big)}\quad\Longrightarrow
\quad n^2(n^2\!-\!P)+n^2m^2\Theta_m(t)^2\le\frac{4\pi^2}{T(E)^2}\qquad\forall t\ge0\, .$$
\end{lemma}
{\em Proof.} From \eq{Econ} and \eq{newTE} we infer that the right hand side of the implication is equivalent to
$$
\Big(n^2-m^2+\sqrt{X_m(E)}\Big)\frac{K\big(\sqrt{Z_m(E)}\big)^2}{\sqrt{X_m(E)}}\le\frac{\pi^2m^2}{4n^2}\, .
$$
By using \eq{firstcase}, we know that $n^2-m^2=(n^2-P)-(m^2-P)\le-(1-\eps^2)(m^2-P)$; therefore, the previous
inequality is certainly satisfied if
$$
\Big(\eps^2+2(1-\eps^2)Z_m(E)\Big)K\big(\sqrt{Z_m(E)}\big)^2\le\frac{\pi^2}{4\eps^2}\, .
$$
This proves the statement.\endproof

For all $0<\eps<1$ and $0<z<\frac12$ we define the function
\neweq{heps}
h_\eps(z):=\max\Bigg\{\frac{\pi^4}{16\eps^4\Big(\eps^2+2(1\!-\!\eps^2)z\Big)^2}\ ,\
\frac{4\, \sigma^4}{\eps^4\, \Big(4(3\eps^4-4\eps^2+\tfrac53)z^2-4(3\eps^4-2\eps^2+\tfrac13)z+3\eps^4\Big)}\Bigg\}
\endeq
and we prove the following bound.

\begin{lemma}\label{bjorn3}
Let $\eps=\tfrac{21}{22}$ and assume that \eqref{firstcase} holds. For all $E>0$ we have
$$K\Big(\sqrt{Z_m(E)}\Big)^4<h_\eps\big(Z_m(E)\big)\, .$$
\end{lemma}
{\em Proof.} We first observe that, for all $\alpha\in(0,\frac{\pi}{2})$,
$$\mbox{the map }z\mapsto\frac{1}{\sqrt{1-z\sin^2\alpha}}\mbox{ is convex for }0<z<\frac12\, .$$
Therefore, the map $f(z):=K(\sqrt{z})$ is convex and, taking into account that $f(0)=\frac{\pi}{2}$ and $f(\frac12 )=\sqrt2\, \sigma$, we infer that
$$
K(\sqrt{z})=f(z)\le\frac{\pi}{2}+(2\sqrt2\, \sigma-\pi)z\qquad\mbox{for }0<z<\frac12\, .
$$
In turn, by taking the fourth power, we obtain that
\neweq{thisone}
K(\sqrt{z})^4\le\left(\frac{\pi}{2}+(2\sqrt2\, \sigma-\pi)z\right)^4\qquad\mbox{for }0<z<\frac12\, .
\endeq
A tedious computation (only involving polynomials) shows that
$$
\left(\frac{\pi}{2}+(2\sqrt2\, \sigma-\pi)z\right)^4<h_\eps(z)\qquad\mbox{for }0<z<\frac12
$$
which, combined with \eq{thisone}, proves the statement.\endproof

Lemma \ref{bjorn3} states that, for all $E>0$, at least one of the implications of Lemmas \ref{bjorn1} and \ref{bjorn} holds.
If the implication of Lemma \ref{bjorn1} holds, then Proposition \ref{lyapzhu} ($ii$) ensures that the trivial solution of
\eq{hill2} is stable; therefore $\Theta_m$ is linearly stable with respect to $\Theta_n$. If the implication of Lemma \ref{bjorn}
holds, then Proposition \ref{lyapzhu} ($i$) leads to the same conclusion. Hence, for all $E>0$, $\Theta_m$ is linearly stable
with respect to $\Theta_n$. This completes the proof of Theorem \ref{stability11}.

\begin{remark}\label{better}
{\em Consider the functions
$$f(z):=K(\sqrt{z})^4\, ,\quad g(z):=\left(\frac{\pi}{2}+(2\sqrt2\, \sigma-\pi)z\right)^4\qquad\mbox{for }0<z<\frac12$$
and, for all $\eps<1$, consider the function $h_\eps$ defined in \eq{heps}.
We proved Lemma \ref{bjorn3} by showing that $f(z)<g(z)<h_\eps(z)$ when $\eps=\frac{21}{22}$. But the function $g$ is not strictly necessary because Lemma \ref{bjorn3} remains true for any $\eps<1$ such that $f(z)<h_\eps(z)$. In Figure \ref{plots} we plot the functions $f$, $g$, and $h_\eps$ for $\eps=\frac{21}{22}$ (left) and $\eps=\frac{26}{27}$ (right).
\begin{figure}[ht]
\begin{center}
\includegraphics[height=4cm,width=8cm]{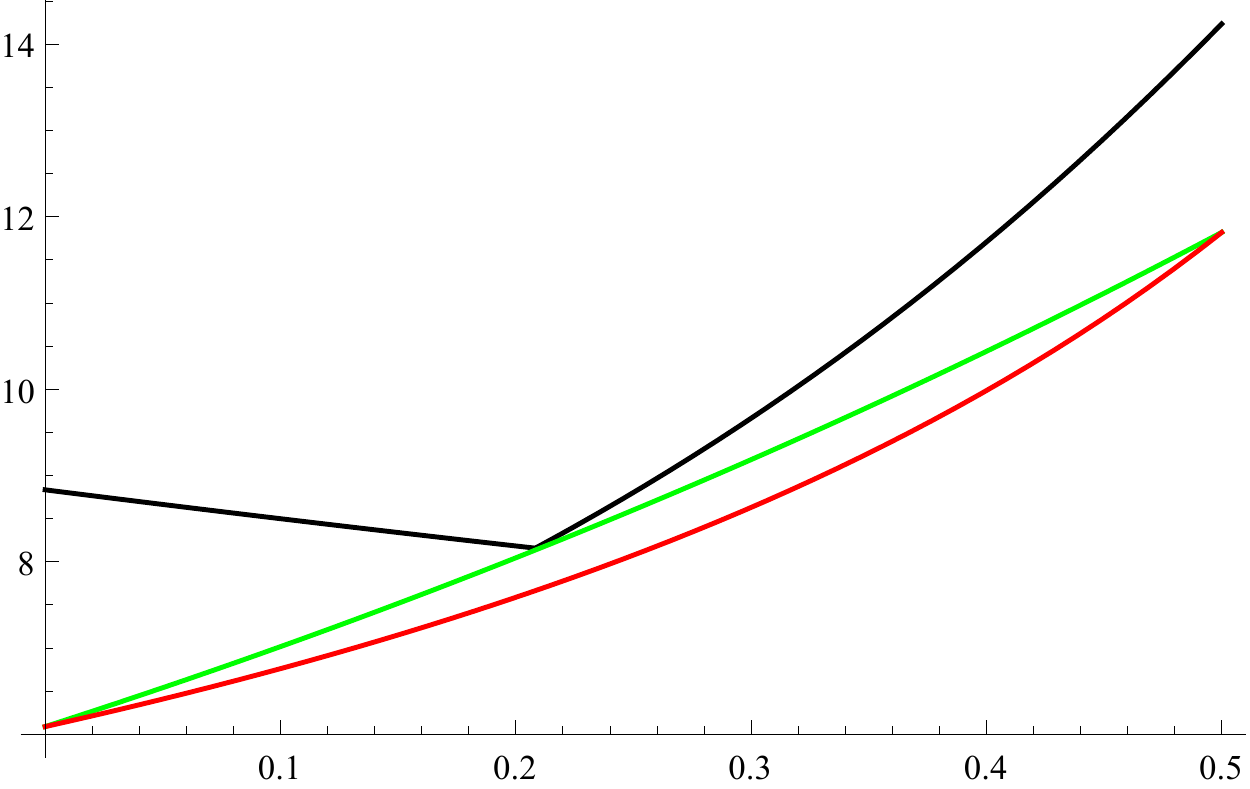}\qquad\includegraphics[height=4cm,width=8cm]{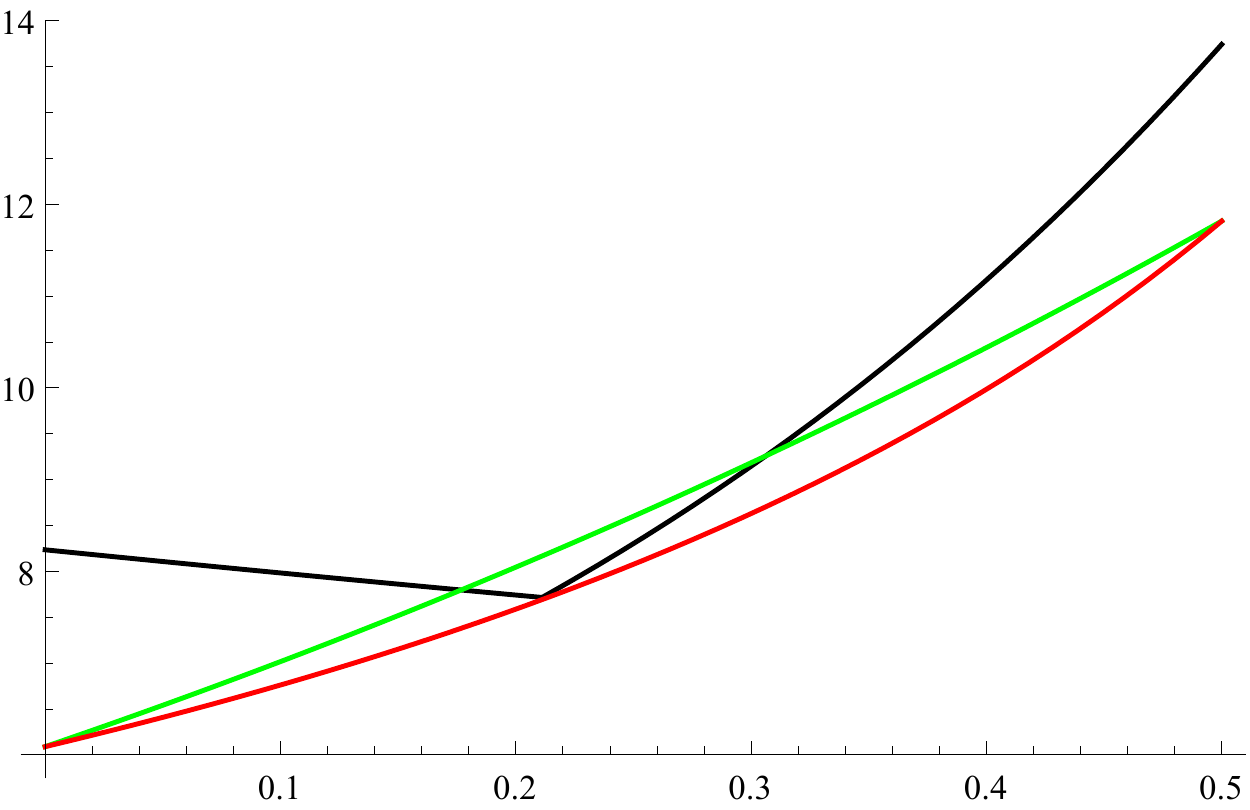}
\caption{Plots of the functions $f$ (red), $g$ (green), and $h_\eps$ (black) for $\eps=\frac{21}{22}$ (left) and $\eps=\frac{26}{27}$ (right).}\label{plots}
\end{center}
\end{figure}
 If we accept Figure \ref{plots} as a proof, then for
$\eps=\frac{26}{27}$ we see that $g(z)<h_\eps(z)$ is no longer true but we still have $f(z)<h_\eps(z)$: then Corollary \ref{conj} may be improved with the bound
$m\le27$. This is the best we can expect from our proof since for $\eps=\frac{27}{28}$ also the inequality $f(z)<h_\eps(z)$ fails for some
$z\in(0,\tfrac12 )$.\par
Finally, note that different theoretical bounds, other than $f(z)<g(z)$ may be obtained by using suitable properties of the complete elliptic integral
of the first kind.\endproof
}\end{remark}
 \section{Proof of Theorems \ref{stability12} and \ref{t:limit-case}}

In this proof we need the following elementary and technical statement:

\begin{lemma}\label{tech1}
Assume that $n^2>m^2>P\ge0$. If there exists an integer $\ell$ such that
\neweq{ell}
n\sqrt{n^2-P}=\ell m\sqrt{m^2-P}\, ,
\endeq
then $P>4m^2-3n^2$.
\end{lemma}
{\em Proof.} From the strict monotonicity of the map $s\mapsto s\sqrt{s^2-P}$ we infer that $\ell>1$, that is, $\ell\ge2$.\par
If it were $n\ge\ell m$, then \eq{ell} would lead to $n\le m$ which contradicts $n>m$. Whence, $n<\ell m$, a fact that we use in the next argument.\par
For contradiction, assume that $P\le4m^2-3n^2$ so that, in particular, $4m^2\ge3n^2$: then
$$
\mbox{\eq{ell}}\Longrightarrow \ell^2m^4\!-\!n^4=P(\ell^2m^2\!-\!n^2)\le(4m^2\!-\!3n^2)(\ell^2m^2\!-\!n^2)
\Longrightarrow 4n^2\ge3\ell^2m^2\ge\frac94 \ell^2n^2\Longrightarrow\frac{16}{9}\ge\ell^2\, ,
$$
which contradicts $\ell\ge2$.\endproof

The existence of an integer $\ell$ as in \eq{ell} is a ``infrequent'' event which, however, may occur: for instance, if $m=1$, $n=2$, $P=3/7$, then $\ell=5$.
As we shall see in the next lemma, this infrequent event deserves a particular attention.\par
The energy estimates of Section \ref{tech} enable us to prove the following statement.

\begin{lemma}\label{premessa}
Assume that $n^2>m^2>P\ge0$. Let $\mu$ be the largest nonnegative
integer such that $\mu m\sqrt{m^2-P}<n\sqrt{n^2-P}$. Then, there
exists $\overline{E}>0$ such that the inequalities
\neweq{doppie}
\mu^2\left(\frac{2\pi}{T(E)}\right)^2\le n^2(n^2-P)+n^2m^2\Theta_m(t)^2\le(\mu+1)^2\left(\frac{2\pi}{T(E)}\right)^2\quad \forall t\ge0
\endeq
are true whenever $0\le E<\overline{E}$.
\end{lemma}
{\em Proof.} By \eq{limitperiod}, when $E=0$ the inequalities in \eq{doppie} become
\neweq{Euguale0}
\mu^2 m^2(m^2-P)\le n^2(n^2-P)\le(\mu+1)^2 m^2(m^2-P)
\endeq
and are therefore fulfilled with {\em strict} inequality on the left. Whence, by continuity, the left inequality in \eq{doppie} remains true for sufficiently small
$E>0$. If also the right inequality in \eq{Euguale0} is strict, then both inequalities in \eq{doppie} remain true for sufficiently small $E$.\par
The only case which remains to be considered is when one has equality on the right of \eq{Euguale0}.
In this case, by Lemma \ref{stimaw} one has for all $t\ge0$:
$$n^2(n^2-P)+n^2m^2\Theta_m(t)^2\le n^2(n^2-P)+n^2m^2\|\Theta_m\|_\infty^2=n^2(n^2-P)+\frac{2 n^2}{m^2-P}\, E+o(E)\qquad\mbox{as }E\to0\, .$$
On the other hand, by \eq{stimavera}, we know that
$$
(\mu+1)^2\left(\frac{2\pi}{T(E)}\right)^2=(\mu+1)^2\left[m^2(m^2-P)+\frac32 \frac{m^2}{m^2-P}\, E\right]+o(E)\, ,\qquad\mbox{as }E\to0
$$
and the statement will follow if we show that
\neweq{claim}
\frac{2 n^2}{m^2-P}<\frac32\, (\mu+1)^2\, \frac{m^2}{m^2-P}=\frac32\, \frac{n^2(n^2-P)}{m^2(m^2-P)}\, \frac{m^2}{m^2-P} \, ,
\endeq
since we assumed that the right inequality in \eq{Euguale0} is an equality. But \eq{claim} is equivalent to
$P>4m^2-3n^2$ which we know to be true after applying Lemma \ref{tech1} with $\ell=\mu+1$. This completes the proof.\endproof

Consider the Hill equation \eq{hill2}: by Theorem \ref{periodicT1}, $a$ is a positive $T(E)/2$-periodic function and by Lemma \ref{premessa} there exists an integer $\mu$ and $\overline{E}>0$ such that
$$
\mu^2\left(\frac{2\pi}{T(E)}\right)^2\le a(t)\le(\mu+1)^2\left(\frac{2\pi}{T(E)}\right)^2\qquad\forall t\ge0
$$
as long as $E\le\overline{E}$. The first criterion in Proposition
\ref{lyapzhu} then states that the trivial solution of \eq{hill2}
is stable. This proves the linear stability for small energies $E$, as stated in Theorem \ref{stability12}-(i).

The statements for large energies of Theorems \ref{stability12} and \ref{t:limit-case} follow from Theorem \ref{stability121}.

\section{Proof of Theorems \ref{stability2} and \ref{stability3}}\label{23}

Assume that $n^2\le \min\{P,m^2\}$ and consider the Hill equation \eq{hill2} where $\Theta_m$ solves \eq{upto} with
$\dot{\Theta}_m(0)=0$. If
\neweq{w0small}
\|\Theta_m\|_\infty^2=\Theta_m(0)^2\le\frac{P-n^2}{m^2}
\endeq
then $n^2(n^2-P)+n^2m^2\Theta_m(t)^2\le0$ for all $t$ and Proposition \ref{negative} states that the trivial solution of
\eq{hill2} is unstable. By \eq{Econ}, the upper bound \eq{w0small} is equivalent to
\neweq{altrastimaE}
E\le-\frac{(P-m^2)^2}4+\frac{(m^2-n^2)^2}4=\frac{(P-n^2)(2m^2-n^2-P)}{4}=:E_1\, .
\endeq

If $n^2<P\le m^2$ (Theorem \ref{stability2}) one has $E_1>0$ while if $n^2<m^2<P$ (Theorem \ref{stability3}) $E_1$ has
the sign of $2m^2-n^2-P$. In any case, when \eq{altrastimaE} holds the trivial solution of \eq{hill2} is unstable and, consequently,
$\Theta_m$ is linearly unstable with respect to $\Theta_n$.\par
Since $m>n$, \eq{this} proves that
$$
\left(\frac{T(E)}{2}\right)^3I(E)<\frac{64}{3}\, \sigma^4
$$
for sufficiently large $E$. Moreover, for large $E$ we also have
$$\int_0^{T(E)/2}\Big(n^2(n^2-P)+n^2m^2 \Theta_m(t)^2\Big)\, dt>0\, .$$
Then the second criterion in Proposition \ref{lyapzhu} ensures
that the trivial solution of \eq{hill2} is stable, provided that
$E$ is sufficiently large, say for $E>E_2$. This proves the linear
stability for $E>E_2$ and completes the proofs of Theorems
\ref{stability2} and \ref{stability3}.

\section{Proof of Theorems \ref{stability22} and \ref{stability4}}

For both theorems, the statements for large energies (ii)-(iii) follow from Theorem \ref{stability121}. Let us prove statement
(i) of both theorems when \eqref{PPP} holds true. If $E_{\Theta_m}<0$, then $0<\sqrt{\Phi_2}\le\Theta_m(t)\le\sqrt{\Phi_1}$ for
all $t$, see Section \ref{pT1}, and the function $a$ in \eq{hill2} is $T(E_{\Theta_m})$-periodic with $T(E_{\Theta_m})$ as given in \eq{TE2}. Moreover,
$$
a(t)\to n^2(n^2-m^2)\qquad\mbox{uniformly as }\ E\downarrow-\frac{(P-m^2)^2}{4}
$$
and, by \eq{limitperiod2},
$$
\lim_{E\downarrow-\frac{(P-m^2)^2}{4}}\left(\frac{\pi}{T(E)}\right)^2=\frac{m^2(P-m^2)}{2}\, .
$$
By \eq{PPP}, there exists an integer $\ell\ge1$ such that
$$\frac{m^2(P-m^2)}{2}\, \ell^2<n^2(n^2-m^2)<\frac{m^2(P-m^2)}{2}\, (\ell+1)^2$$
and, by continuity, there exists $-\frac{(P-m^2)^2}{4}<E_1\le0$ such that
$$\left(\frac{\pi}{T(E)}\right)^2\, \ell^2\le a(t)\le\left(\frac{\pi}{T(E)}\right)^2\, (\ell+1)^2\qquad\forall t\ge0$$
whenever $-\frac{(P-m^2)^2}{4}<E_{\Theta_m}<E_1$. Then the first
criterion in Proposition \ref{lyapzhu} ensures that the trivial
solution of \eq{hill2} is stable if $E_{\Theta_m}$ belongs to this interval.

Next we turn to the much more involved proof of (i) for both theorems when \eqref{PPP-2}
holds true. We divide the proof in several steps.\par \noindent
\emph{Step 1: asymptotic behavior of the solution of $\eq{ODE}$ for negative energies.} Let us define $\overline E:=-\frac{(P-m^2)^2}4$ and
let us put $\eps=\sqrt{E-\overline E}$ with $E\ge \overline E$. We
denote by $u_\eps$ the solution of the Cauchy problem
\begin{equation}\label{eq:u-eps}
\begin{cases}
\ddot{u}_\eps+m^2(m^2-P)u_\eps+m^4 u_\eps^3=0 \, , \\[8pt]
u_\eps(0)=\sqrt{\frac{P-m^2}{m^2}+\frac 2{m^2}\eps} \, , \qquad
\dot{u}_\eps(0)=0 \, .
\end{cases}
\end{equation}
We observe that $u_\eps$ is a solution at energy $E=\overline{E}+\eps^2$ and, by using standard arguments from the theory of ordinary
differential equations, the map $\eps \mapsto u_\eps(t)$ is
smooth. We claim that
\begin{equation}  \label{eq:Taylor-u-eps-quadro}
u_\eps(t)^2=\tfrac{P-m^2}{m^2}+\tfrac{2\cos(\sqrt 2
m\sqrt{P-m^2}t)}{m^2}\, \eps -\tfrac{2\sin^2(\sqrt 2
m\sqrt{P-m^2}t)}{m^2(P-m^2)} \, \eps^2+o(\eps^2) \qquad \text{as }
\eps\to 0^+ \, ,
\end{equation}
uniformly on bounded time intervals. In order to show this, we use the notation
$u(t,\eps):=u_\eps(t)$, $A_\eps(t):=\partial_\eps u(t,\eps)$ and
$B_\eps(t):=\partial^2_\eps u(t,\eps)$. The functions $A_\eps$ and $B_\eps$ for $\eps=0$ can be obtained
explicitly by solving the following linear Cauchy problems coming
from formal differentiations with respect to $\eps$ in \eqref{eq:u-eps}:
\begin{align*}
\begin{cases}
\ddot{A}_0(t)+2m^2(P-m^2)A_0(t)=0 \, ,\\[7pt]
A_0(0)=\frac 1{m\sqrt{P-m^2}} \, , \qquad \dot{A}_0(0)=0
\end{cases}
\end{align*}
and
\begin{align*}
\begin{cases}
\ddot{B}_0(t)+2m^2(P-m^2)B_0(t)+\frac{6m}{\sqrt{P-m^2}} \, \cos^2(\sqrt 2 m \sqrt{P-m^2}t)=0 \, , \\[7pt]
B_0(0)=-\frac 1{m(P-m^2)^{3/2}} \, , \qquad \dot{B}_0(0)=0 \, .
\end{cases}
\end{align*}
From $A_0$ and $B_0$ we readily obtain the second order Taylor expansion as $\eps\to 0^+$ of $u_\eps$:
\begin{equation*}
u_\eps(t)=\sqrt{\tfrac{P-m^2}{m^2}}+\tfrac{\cos(\sqrt 2
m\sqrt{P-m^2}t)}{m\sqrt{P-m^2}}\, \eps +\tfrac{\cos(2\sqrt 2
m\sqrt{P-m^2}t)-3}{4m(P-m^2)^{3/2}}\, \eps^2 +o(\eps^2) \qquad
\text{as } \eps\to 0^+ \, ,
\end{equation*}
uniformly on bounded time intervals. By squaring we reach \eq{eq:Taylor-u-eps-quadro}.
\par \noindent
\emph{Step 2: switch to polar coordinates.} For any $\eps>0$ and $x,y\in \R$, we define $v_\eps$ as the
unique solution of the Cauchy problem
\begin{equation} \label{eq:v-eps}
\begin{cases}
\ddot{v}_\eps+n^2(n^2-P)v_\eps+n^2m^2 u_\eps^2 v_\eps=0 \, , \\
v_\eps(0)=x \, , \qquad \dot{v}_\eps(0)=\sqrt{n^2(n^2-m^2)}\, y \,
.
\end{cases}
\end{equation}
We put
$x_\eps(t):=v_\eps(t)$ and $y_\eps(t):=(n^2(n^2-m^2))^{-1/2}
\dot{x}_\eps(t)$. In order to better understand the behavior of the
solution in the $xy$-plane, we switch to polar coordinates by
defining     the functions $\rho_\eps$, $\theta_\eps$ in such a
way that $x_\eps(t)=\rho_\eps(t)\cos(\theta_\eps(t))$ and
$y_\eps(t)=\rho_\eps(t)\sin(\theta_\eps(t))$. Then, by \eq{eq:v-eps} we obtain
\begin{equation} \label{eq:rho-theta-eps}
\begin{cases}
\dot{\rho}_\eps=\frac{n^2(P-m^2)-n^2m^2 u_\eps^2}{\sqrt{n^2(n^2-m^2)}} \, \rho_\eps\sin(\theta_\eps)\cos(\theta_\eps) \, , \\[10pt]
\dot \theta_\eps=-\sqrt{n^2(n^2-m^2)}\sin^2(\theta_\eps)-\frac{n^2(n^2-P)+n^2m^2u_\eps^2}{\sqrt{n^2(n^2-m^2)}}\cos^2(\theta_\eps)\, , \\
\rho_\eps(0)=\rho \, , \qquad \theta_\eps(0)=\theta \, ,
\end{cases}
\end{equation}
where $\rho,\theta$ are such that $x=\rho\cos\theta$ and
$y=\rho\sin\theta$.

Then we introduce the functions
$\rho_0,\rho_1,\rho_2,\theta_0,\theta_1,\theta_2$ such that
\begin{align*}
\rho_\eps(t)=\rho_0(t)+\rho_1(t)\eps+\rho_2(t)\eps^2+o(\eps^2) \,
, \quad
\theta_\eps(t)=\theta_0(t)+\theta_1(t)\eps+\theta_2(t)\eps^2+o(\eps^2)
\, , \quad \text{as } \eps \to 0^+
\end{align*}
uniformly on bounded time intervals. The existence of such
functions can be proved by showing that $\rho_\eps$ and
$\theta_\eps$ are smooth with respect to $\eps$ at
$\eps=0$.
\par \noindent
\emph{Step 3: characterization of the functions $\rho_0,\rho_1,\rho_2,\theta_0,\theta_1,\theta_2$.} To compute explicitly $\rho_0$ and $\theta_0$, it is
sufficient to choose $\eps=0$ in \eqref{eq:rho-theta-eps} and to exploit the fact that $u_0\equiv\sqrt{\frac{P-m^2}{m^2}}$. Then one obtains
\begin{align*}
\dot{\rho}_0=0 \quad \text{and} \quad \dot{\theta}_0=-\sqrt{n^2(n^2-m^2)}\,,
\end{align*}
from which it follows that
\begin{align} \label{eq:rho0-theta0}
\rho_0(t)=\rho \quad \text{and} \quad
\theta_0(t)=\theta-\sqrt{n^2(n^2-m^2)} \, t \qquad \text{for any }
t \, .
\end{align}

Similarly, in order to compute $\rho_1$ and $\theta_1$, one has to
differentiate with respect to $\eps$ in \eqref{eq:rho-theta-eps}
and to put $\eps=0$. Taking into account \eqref{eq:Taylor-u-eps-quadro}
we obtain
\begin{align*}
& \dot{\rho}_1=-\frac{2n}{\sqrt{n^2-m^2}} \, \cos(\sqrt 2 m \sqrt{P-m^2}\, t) \rho_0 \sin(\theta_0)\cos(\theta_0) \, , \\
& \dot{\theta}_1=-\frac{2n}{\sqrt{n^2-m^2}} \, \cos(\sqrt 2 m
\sqrt{P-m^2}\, t)\cos^2(\theta_0) \, .
\end{align*}
Therefore, since $\rho_1(0)=0$ and $\theta_1(0)=0$, by
\eqref{eq:rho0-theta0} we obtain
\begin{align*}
& \rho_1(t)=-\tfrac{2n\rho}{\sqrt{n^2-m^2}} \int_0^t  \cos\big(\sqrt 2 m \sqrt{P-m^2}\, s\big) \sin\big(\theta-\sqrt{n^2(n^2-m^2)} \, s\big) \cos\big(\theta-\sqrt{n^2(n^2-m^2)} \, s\big)\, ds\,, \\
& \theta_1(t)=-\tfrac{2n}{\sqrt{n^2-m^2}} \int_0^t  \cos\big(\sqrt
2 m \sqrt{P-m^2}\, s\big)\cos^2\big(\theta-\sqrt{n^2(n^2-m^2)} \,
s\big)\, ds\,.
\end{align*}
To simplify the notation we define
\begin{equation} \label{eq:def-ABC}
A=A(P,m):=\sqrt 2 m \sqrt{P-m^2} \, , \quad  B:=\theta \, , \quad
C=C(n,m):=\sqrt{n^2(n^2-m^2)}\, .
\end{equation}
After some computations one obtains
\begin{align} \label{eq:rho1-theta1}
& \rho_1(t)=-\tfrac{2n\rho}{\sqrt{n^2-m^2}}
\left\{\tfrac{\cos\big[(A+2C)t-2B\big]-\cos(2B)}{4(A+2C)}
-\tfrac{\cos\big[(A-2C)t+2B\big]-\cos(2B)}{4(A-2C)}\right\} \, , \\
\notag &
\theta_1(t)=-\tfrac{2n}{\sqrt{n^2-m^2}}\left\{\tfrac{\sin(At)}{2A}+\tfrac{\sin\big[(A+2C)t-2B\big]+\sin(2B)}{4(A+2C)}
+\tfrac{\sin\big[(A-2C)t+2B\big]-\sin(2B)}{4(A-2C)}\right\} \, .
\end{align}

Finally, we recover an explicit representation for
$\rho_2$ and $\theta_2$. We start by expanding the following term which
appears in the first equation in \eqref{eq:rho-theta-eps}
\begin{align*}
\rho_\eps(t)\sin(\theta_\eps(t))\cos(\theta_\eps(t))=\tfrac{\rho_0(t)}2\sin(2\theta_0(t))
+\left[\tfrac{\rho_1(t)}2 \sin(2\theta_0(t))+\rho_0(t)
\cos(2\theta_0(t))\theta_1(t)\right]\eps+o(\eps)
\end{align*}
as $\eps\to 0$ uniformly on bounded time intervals. By \eqref{eq:Taylor-u-eps-quadro} and \eqref{eq:rho0-theta0}-\eqref{eq:rho1-theta1} we obtain
\begin{align*}
& \dot{\rho}_2=-\tfrac{2n}{\sqrt{n^2-m^2}}\,
\cos(At)\left[\tfrac{\rho_1}2 \sin(2\theta_0)+\rho_0
\cos(2\theta_0)\theta_1\right]
+\tfrac{n\sin^2(At)}{(P-m^2)\sqrt{n^2-m^2}} \, \rho_0\sin(2\theta_0)\,, \\[7pt]
& \dot{\theta}_2=\tfrac{2n}{\sqrt{n^2-m^2}}\,
\cos(At)\sin(2\theta_0)\theta_1+\tfrac{2n}{(P-m^2)\sqrt{n^2-m^2}}
\, \sin^2(At) \cos^2(\theta_0)\,,
\end{align*}
so that
\begin{align} \label{eq:rho2-theta2}
& \rho_2(t)=\int_0^t \left\{-\tfrac{2n}{\sqrt{n^2-m^2}}\,
\cos(As)\left[\tfrac{\rho_1}2 \sin(2\theta_0)+\rho_0\cos(2\theta_0)\theta_1\right]
+\tfrac{n\sin^2(As)}{(P-m^2)\sqrt{n^2-m^2}} \, \rho_0\sin(2\theta_0)\right\} \, ds \, , \\
\notag & \theta_2(t)=\int_0^t \left\{\tfrac{2n}{\sqrt{n^2-m^2}}\,
\cos(As)\sin(2\theta_0)\theta_1+\tfrac{2n}{(P-m^2)\sqrt{n^2-m^2}}
\, \sin^2(As) \cos^2(\theta_0) \right\} \, ds \, .
\end{align}
We observe that the two above integrals can be computed by using the explicit representations of
$\rho_0,\rho_1,\theta_0,\theta_1$ given in \eqref{eq:rho0-theta0} and \eqref{eq:rho1-theta1} but, as we will see below, we only need
to compute their values at $T_\eps:=T(\overline E+\eps^2)$ where $T(E)$ is given by \eqref{TE2}.
\par \noindent
\emph{Step 4: asymptotic behavior of $T_\eps$.}
By \eqref{TE2} we have
\begin{equation*}
T_\eps=\frac{2\sqrt 2}{m\sqrt{P-m^2+2\eps}} \int_{\delta_\eps}^1
\frac{1}{\sqrt{(1-s^2)(s^2-\delta_\eps^2)}}\, ds
\end{equation*}
where we put $\delta_\eps:=\left(1-\frac2{P-m^2}\eps\right)\left(1-\frac 4{(P-m^2)^2} \, \eps^2\right)^{-1/2}$. We first observe that
\begin{equation} \label{eq:T-1}
\frac{2\sqrt 2}{m\sqrt{P-m^2+2\eps}}=\frac{2\sqrt
2}{m\sqrt{P-m^2}}-\frac{2\sqrt 2}{m(P-m^2)^{3/2}}\,
\eps+\frac{3\sqrt 2}{m(P-m^2)^{5/2}}\, \eps^2 +o(\eps^2) \, ,
\qquad \text{as } \eps\to 0^+ \, .
\end{equation}
Moreover, we also have
\begin{equation*}
\delta_\eps=1-\frac 2{P-m^2} \, \eps+\frac 2{(P-m^2)^2}\,
\eps^2+o(\eps^2) \, , \qquad \text{as } \eps\to 0^+ \, .
\end{equation*}

If we define $f(\delta):=\int_\delta^1  \frac
1{\sqrt{(1-s^2)(s^2-\delta^2)}} \, ds$ for any $\delta\in (0,1)$,
we have that
\begin{equation} \label{eq:T-3}
f(\delta)=\frac \pi 2+\frac \pi
4(1-\delta)+\frac{5\pi}{32}(1-\delta)^2+o((1-\delta)^2) \qquad
\text{as } \delta\to 1^- \, .
\end{equation}
To see this, we compute $f(1):=\lim_{\delta\to 1^-}f(\delta)$, $f'(1)$
and $f''(1)$, by first writing $f$ in the form
\begin{equation*}
f(\delta)=\frac{2}{1-\delta} \int_\delta^1 \frac 1{\sqrt{1-\left[\frac 2{1-\delta}\left(s-\frac{\delta+1}2\right)\right]^2}} \cdot \frac 1{\sqrt{(1+s)(s+\delta)}} \, ds
\end{equation*}
and then using the change of variable $t=\frac 2{1-\delta}\left(s-\frac{\delta+1}2\right)$ to get
\begin{equation*}
f(\delta)=\int_{-1}^1 \frac 1{\sqrt{1-t^2}}
\frac{1}{\sqrt{\left(\frac{1-\delta}2\,
t+\frac{\delta+3}2\right)\left(\frac{1-\delta}2\,
t+\frac{3\delta+1}2\right)}}\, dt \, .
\end{equation*}
Combining \eqref{eq:T-1}-\eqref{eq:T-3} we obtain
\begin{equation} \label{eq:T-5}
T_\eps=\frac{\pi\sqrt 2}{m\sqrt{P-m^2}}+\frac{3\pi \sqrt
2}{4m(P-m^2)^{5/2}}\, \eps^2+o(\eps^2)  \qquad \text{as } \eps\to
0^+ \, .
\end{equation}
\par \noindent
\emph{Step 5: evaluation of the functions $\rho_0,\rho_1,\rho_2,\theta_0,\theta_1,\theta_2$ at $T_\eps$.} Taking into account that the first order term in the asymptotic
expansion of $T_\eps$ vanishes, see \eq{eq:T-5}, by direct computation one sees that
\begin{equation} \label{eq:rho1-theta1-T/2}
\begin{split}
& \rho_1(T_\eps)=\rho_1(T_0)+o(\eps)=\tfrac{n^2 \rho}{m^2(P-m^2)(1-L^2)}
\left[\cos(2\pi L-2\theta)-\cos(2\theta)\right]+o(\eps) \quad \text{as } \eps\to 0^+ \, ,\\[5pt]
 & \theta_1(T_\eps)=\theta_1(T_0)+o(\eps)=\tfrac{n^2}{m^2(P-m^2)(1-L^2)}[\sin(2\pi
L-2\theta)+\sin(2\theta)]+o(\eps) \quad \text{as } \eps\to 0^+ \, ,
\end{split}
\end{equation}
where $L$ is the number defined in \eqref{PPP} and $T_0:=\lim_{\eps\to 0^+} T_\eps$. By \eqref{eq:T-5} and \eqref{eq:rho1-theta1-T/2}, we may write
\begin{align} \label{eq:expansions}
\begin{split}
& \rho_\eps(T_\eps)=\rho+\rho_1(T_0)\eps+\rho_2(T_0)\eps^2+o(\eps^2) \quad \text{as } \eps\to 0^+ \, ,\\[7pt]
 & \theta_\eps(T_\eps)=\theta-\pi L+\theta_1(T_0)\eps+\theta_2(T_0)\eps^2-\tfrac{3\pi
L}{4(P-m^2)^2}\, \eps^2+o(\eps^2)\quad \text{as } \eps\to 0^+ \, .
\end{split}
\end{align}
In turn, by \eqref{eq:rho0-theta0} and \eqref{eq:expansions} we obtain
\begin{align} \label{eq:v-eps(T-eps)}
 x_\eps(T_\eps)&=\rho_\eps(T_\eps)\cos(\theta_\eps(T_\eps))\\
\notag &
=\left[\rho+\rho_1(T_0)\eps+\rho_2(T_0)\eps^2+o(\eps^2)\right]
\Big\{\cos(\theta-\pi L)-\sin(\theta-\pi L)\theta_1(T_0)\eps \\
\notag &\quad \left. +\left[\sin(\theta-\pi L)\left(\tfrac{3\pi
L}{4(P-m^2)^2}-\theta_2(T_0)\right) -\tfrac 12 \cos(\theta-\pi
L)(\theta_1(T_0))^2 \right]\eps^2+o(\eps^2) \right\}
\end{align}
and
\begin{align} \label{eq:v'-eps(T-eps)}
 y_\eps(T_\eps)&=\rho_\eps(T_\eps)\sin(\theta_\eps(T_\eps))\\
\notag &
=\left[\rho+\rho_1(T_0)\eps+\rho_2(T_0)\eps^2+o(\eps^2)\right]
\Big\{\sin(\theta-\pi L)+\cos(\theta-\pi L)\theta_1(T_0)\eps \\
\notag & \quad \left. +\left[\cos(\theta-\pi
L)\left(\theta_2(T_0)-\tfrac{3\pi L}{4(P-m^2)^2}\right) -\tfrac 12
\sin(\theta-\pi L)(\theta_1(T_0))^2 \right]\eps^2+o(\eps^2)
\right\}\,.
\end{align}
\par \noindent
\emph{Step 6: introduction of the monodromy matrix.}
Let us denote by $M_\eps$ the monodromy matrix of \eq{eq:v-eps}, see \cite[Chapter
II, Section 2.1]{yakubovich} for its precise definition. In our case we have
\begin{equation*}
M_\eps=
\begin{pmatrix}
x_\eps(T_\eps)_{|(\rho=1, \theta=0)} & x_\eps(T_\eps)_{|(\rho=1, \theta=\pi/2)} \\
y_\eps(T_\eps)_{|(\rho=1, \theta=0)} & y_\eps(T_\eps)_{|(\rho=1,
\theta=\pi/2)}
\end{pmatrix}\ .
\end{equation*}
By inserting \eqref{eq:v-eps(T-eps)}-\eqref{eq:v'-eps(T-eps)} into $M_\eps$,
for any $L>1$ not necessarily integer, we infer that
\begin{equation*}
M_\eps=
\begin{pmatrix}
\cos(\pi L)+o(1)  & \sin(\pi L)+o(1) \\
-\sin(\pi L)+o(1) &  \cos(\pi L)+o(1)
\end{pmatrix} \ .
\end{equation*}
If $L\not\in \N$ and $\eps$ is small enough, then the
eigenvalues of $M_\eps$ are complex numbers with nontrivial
imaginary part, thus we recover statement (i) of both Theorems \ref{stability22} and \ref{stability4} when \eqref{PPP} holds true.
\par \noindent
\emph{Step 7: asymptotic behavior of $M_\eps$.}
In the last part of the proof we use \eq{PPP-2} from which we infer that $L>1$ ($L\in \N$). To obtain an expansion for each
component of the matrix $M_\eps$ we may assume that
$\rho=1$ and $\theta=0$ or $\theta=\pi/2$ for, respectively, the
components of the first and the second column of $M_\eps$. By
\eqref{eq:rho1-theta1-T/2} one sees that
$\rho_1(T_0)=\theta_1(T_0)=0$ both for $\theta=0$ and
$\theta=\pi/2$.

We now claim that $\rho_2(T_0)=0$ both for $\theta=0$ and
$\theta=\pi/2$. To see this, one has to insert
\eqref{eq:rho0-theta0} and \eqref{eq:rho1-theta1} into
\eqref{eq:rho2-theta2} and compute explicitly all the integrals.
Since the computations only involve elementary calculus, we omit them; let us just mention that all these integrals
vanish since, by using Werner formulas, all of them may be reduced to
an integral of the type $\int_0^{T_0} \sin[(k_1 A+2k_2 C)t] \, dt$
both for $\theta=0$ and $\theta=\pi/2$, where $A, C$ are defined
in \eqref{eq:def-ABC} and $k_1, k_2\in \Z$. Finally, since $A
T_0=2\pi$ and $CT_0=\pi L$, then $\int_0^{T_0} \sin[(k_1 A+2k_2 C)t] \, dt=0$ whenever $L\in \N$.

Let us now compute $\theta_2(T_0)$.
Differently from the computation of $\rho_2(T_0)$, not all
the integrals coming from \eqref{eq:rho0-theta0},
\eqref{eq:rho1-theta1}, \eqref{eq:rho2-theta2} vanish, even if
$L\in \N$. Some of them are equal to $\frac{T_0}4$, $-\frac{T_0}4$, and
$\frac{T_0}2$. After some tedious computations, one gets
\begin{align*}
& \theta_2(T_0)=\tfrac{2n}{\sqrt{n^2-m^2}}
\left\{-\tfrac{2n}{\sqrt{n^2-m^2}}\left[-\tfrac{1}{4(A+2C)}\tfrac{T_0}4+\tfrac{1}{4(A-2C)}\tfrac{T_0}4\right]\right\}
+\tfrac{2n}{(P-m^2)\sqrt{n^2-m^2}}\tfrac{T_0}4 \\[8pt]
& \qquad =-\tfrac{n^2 T_0
L}{2(n^2-m^2)A(1-L^2)}+\tfrac{nT_0}{(P-m^2)\sqrt{n^2-m^2}}=\tfrac{2
n^3 T_0}{4\sqrt{n^2-m^2}m^2
(P-m^2)(L^2-1)}+\tfrac{nT_0}{(P-m^2)\sqrt{n^2-m^2}}\,.
\end{align*}
Coming back to \eqref{eq:v-eps(T-eps)} and \eqref{eq:v'-eps(T-eps)}
we obtain
\begin{equation} \label{eq:M-eps-2}
M_\eps=
\begin{pmatrix}
(-1)^L+o(\eps^2)   &  -(-1)^L\left(\theta_2(T_0)-\tfrac{3\pi L}{4(P-m^2)^2}\right)\eps^2+o(\eps^2) \\
(-1)^L\left(\theta_2(T_0)-\tfrac{3\pi
L}{4(P-m^2)^2}\right)\eps^2+o(\eps^2) & (-1)^L+o(\eps^2)
\end{pmatrix}\,.
\end{equation}
\par \noindent
\emph{Step 8: eigenvalues of $M_\eps$ and conclusions.}
Now we observe that if $v$ is a solution of the Hill equation in
\eqref{eq:v-eps} then the function $t\mapsto v(T_\eps-t)$ is a
solution of the same equation. This yields that, for any
$a,b\in\R$, the following implication holds
\begin{equation*}
\left(
\begin{array}{rr}
\tilde a \\ \tilde b
\end{array}\right)
=M_\eps \left(
\begin{array}{rr}
a \\ b
\end{array}
\right) \quad\Longrightarrow \quad
M_\eps \left(
\begin{array}{rr}
\tilde a \\ -\tilde b
\end{array}
\right)=\left(
\begin{array}{rr}
a \\ -b
\end{array}\right)\,.
\end{equation*}
Proceeding similarly to the proof of \cite[Lemma
3.3]{cazw2}, we infer that ${\rm det}(M_\eps)=1$ and the diagonal components of $M_\eps$ are equal, namely $(M_\eps)_{11}=(M_\eps)_{22}$.

We define
$g(\eps):=(M_\eps)_{11}-(-1)^L=(M_\eps)_{22}-(-1)^L$ in such a way
that $g(\eps)=o(\eps^2)$ as $\eps\to 0^+$, as one can deduce from
\eqref{eq:M-eps-2}. We put
$K:=(-1)^L\left(\theta_2(T_0)-\tfrac{3\pi L}{4(P-m^2)^2}\right)$ and we obtain
\begin{equation*}
1={\rm det}(M_\eps)=1+2(-1)^L g(\eps)+K^2 \eps^4+o(\eps^4)
\end{equation*}
and hence
\begin{equation*}
\frac{g(\eps)}{\eps^4}=\tfrac{(-1)^{L+1}}2 K^2+o(1)  \, .
\end{equation*}
\par \noindent
We observe that we may rewrite $K$ in the form
\begin{align*}
K&=(-1)^L \tfrac{nT_0}{4\sqrt{n^2-m^2}(P-m^2)}
\left[\tfrac{2n^2}{m^2(L^2-1)}+2-\tfrac{3m^2}{2n^2}\, L^2\right]
\\
 \qquad &=(-1)^L \left(\tfrac{nT_0}{4\sqrt{n^2-m^2}(P-m^2)} \,
\tfrac{3m^4L^4-(3m^4+4n^2m^2)L^2+4n^2m^2-4n^4}{2n^2m^2(1-L^2)}\right)
\, .
\end{align*}
Therefore, \eqref{PPP-2} implies $K \neq 0$ so that $g(\eps)$ is
eventually negative as $\eps\to 0^+$ when $L$ is even and
eventually positive as $\eps\to 0^+$ when $L$ is odd.

But $|{\rm tr}(M_\eps)|=2+2(-1)^L g(\eps)<2$ eventually as
$\eps\to 0^+$. Recalling again that ${\rm det}(M_\eps)=1$, this
implies that the eigenvalues of $M_\eps$ are necessarily complex conjugate with
nontrivial imaginary part. By the Floquet theory this completes the
proof of (i) of Theorems \ref{stability22} and \ref{stability4}, see e.g.\ \cite[Chapter
II]{yakubovich}.

\section{Appendix: complements and computations}
\subsection{The stationary problem}

Stationary solutions of \eq{truebeam} solve the following boundary value problem:
\neweq{nonlineareq}
\left\{\begin{array}{l}
u''''+\Big[P-\frac{2}{\pi}\, \|u'\|^2_{L^2(0,\pi)}\Big]u''=0\quad x\in(0,\pi)\, ,\\
u(0)=u(\pi)=u''(0)=u''(\pi)=0\, .
\end{array}\right.
\endeq
This problem has been studied by many authors from different point of view: the most relevant contributions related to our paper are \cite{grotta, yaga1}. Solutions of \eq{nonlineareq} are critical points
of the functional
$$J_0(v)=\frac12 \int_0^\pi(v'')^2-\frac{P}{2}\int_0^\pi(v')^2+\frac{1}{2\pi}\left(\int_0^\pi(v')^2\right)^2\qquad(v\in H^2\cap H^1_0(0,\pi))\, .$$
The following precise multiplicity result for \eq{nonlineareq} holds:

\begin{proposition}\label{homogeneous}
If $P\in(k^2,(k+1)^2]$ for some $k\ge0$, then \eqref{nonlineareq} admits exactly $2k+1$ solutions which are explicitly given by
\begin{equation}  \label{onlyeigenfunctions}
u_0(x)=0\,,\qquad\pm u_j(x)=\pm\, \frac{\sqrt{P-j^2}}{j}\ \sin(jx)\quad(j=1,...,k)\, .
\end{equation}
Moreover, for each solution the energy is given by
\begin{equation}  \label{energysolutions}
J_0(u_0)=0\,,\qquad J_0(\pm u_j)=-\frac{\pi}{8}(P-j^2)^2\quad(j=1,...,k)\, ,
\end{equation}
and the Morse index $M$ is given by
\begin{equation}  \label{morse}
M(u_0)=k\,,\qquad M(\pm u_j)=j-1\quad(j=1,...,k)\,.
\end{equation}
\end{proposition}\par\bigskip

\begin{figure}[ht]
\begin{center}
 {\includegraphics[height=25mm, width=80mm]{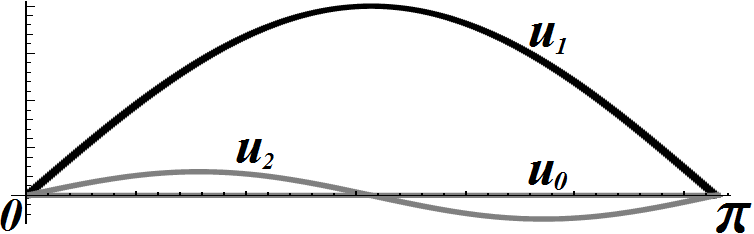}}
\caption{Equilibrium positions of the compressed beam when $P=5$.}\label{twobeams}
\end{center}
\end{figure}

The proof of Proposition \ref{homogeneous} may be obtained by combining results from Reiss \cite{reiss}, see also \cite[Section 3]{reissm}, with results from \cite{algwaiz}. For any $P>1$, formula \eq{morse} in Proposition \ref{homogeneous} states that only $\pm u_1$ are stable (global minima) while all the
other critical points of $J_0$ are saddle points. In particular, if the beam is subject to the compression $P=5$, then the solutions are
$$u_0(x)=0\, ,\qquad\pm u_1(x)=\pm2 \, \sin(x)\, ,\qquad\pm u_2(x)=\pm \frac{1}{2}\ \sin(2x)\, ,$$
see \eq{onlyeigenfunctions}. The solutions $u_0$, $u_1$, and $u_2$ are depicted in scale in Figure \ref{twobeams}.
The minimum of $J_0$ is attained only at $\pm u_1$: from a physical point of view, this means that the compressed beam in position $u_1$ is
stable while the beams in positions $u_0$ and $u_2$ are not.\par

\subsection{Some comments about the bounds for the energy}\label{some}

If $P\le k^2$ then \eq{ODE} admits no solutions with negative energy, while if $E(\alpha,\beta)=0$ then $\alpha=\beta=0$ and the
solution of \eq{ODE}-\eq{alphabeta} is trivial: $\Theta_k(t)\equiv0$.\par If $P>k^2$, then necessarily
$E\ge-\frac{(P-k^2)^2}{4}$ with equality only for the two constant solutions $\Theta_k(t)\equiv\pm\frac{\sqrt{P-k^2}}{k}$. Moreover,
if $E=0$ then \eq{ODE}-\eq{alphabeta} may admit non-constant and non-periodic solutions: the function
$$
\overline{\Theta}_k(t)=\frac{\sqrt2\, \sqrt{P-k^2}}{k\, \cosh(k\sqrt{P-k^2}\, t)}\ ,
$$
as well as the opposite function $-\overline{\Theta}_k(t)$ and their translations $\overline{\Theta}_k(t-t_0)$ and $-\overline{\Theta}_k(t-t_0)$ for
any $t_0\in\R$, solve \eq{ODE}-\eq{alphabeta} for different values of $\alpha$ and $\beta$. They all have energy $E=0$ and are
homoclinic to $0$: they may be observed in the usual $\infty$-shaped picture for the orbits in the phase plane.

\subsection{Stability criteria for the Hill equation}\label{stabcrit}

Throughout this paper we made use of some properties of the Hill equation, that is
\neweq{hill}
\ddot{\xi}+a(t)\xi=0\, ,\qquad a\in C^0([0,T])\, ,\quad a(t+T)=a(t)\quad \forall t
\endeq
where we intend that $T>0$ is the smallest period of $a$. This equation was introduced by Hill \cite{Hill} for the study of the lunar perigee and
has been the object of many subsequent studies, see e.g.\ \cite{cesari,chicone,magnus,stoker}. The main concern is to establish whether the trivial
solution $\xi\equiv0$ of \eq{hill} is stable or, equivalently, if all the solutions of \eq{hill} are bounded in $\R$. The following stability criteria
have been used in the course.

\begin{proposition}\label{lyapzhu}
Let $\sigma$ be as in \eqref{sigma}. Assume that one of the two following facts holds:\par
$$(i)\quad a\ge0\quad\mbox{and}\quad
\exists\ell\in\N\quad\mbox{s.t.}\quad\frac{\ell^2\pi^2}{T^2}\le a(t)\le\frac{(\ell+1)^2\pi^2}{T^2}\quad\forall t\, ,$$
$$(ii)\quad\int_0^Ta(t)\, dt>0\quad\mbox{and}\quad T^3\int_0^Ta^+(t)^2\, dt<\frac{64}{3}\, \sigma^4\, ;$$
then the trivial solution of \eqref{hill} is stable (here, $a^+=\max\{a,0\}$).
\end{proposition}

The first criterion is due to Zhukovskii \cite{zhk} (see also \cite[Chapter VIII]{yakubovich}). The second criterion is due to
Li-Zhang \cite[Theorem 1]{lizhang} (case $\alpha=2$) and generalizes classical criteria by Lyapunov \cite{lyapunov} and
Borg \cite{borg}. The criteria in Proposition \ref{lyapzhu} are somehow ``dual'': $(i)$ is needed for small energies while $(ii)$
is needed for large energies.\par Concerning instability, we state a simple sufficient condition, see (4.2.i) p.60 in \cite{cesari}.

\begin{proposition}\label{negative}
Assume that $a\le0$; then the trivial solution of \eqref{hill} is unstable.
\end{proposition}

Note that Proposition \ref{negative} cannot be relaxed with the requirement that $a$ has negative mean
value ($\int_0^Ta\le0$), see the Corollary on p.697 in \cite{yakubovich}.\par\bigskip\noindent

\textbf{Acknowledgments.} The authors are grateful to Fabio Zanolin for several discussions and to anonymous referees that allowed to improve the first version of this paper.
The first and second author are partially supported by the FIR 2013 project \emph{Geometrical and qualitative aspects of PDE's}.
The third and fourth authors are partially supported by the PRIN project {\em Equazioni alle derivate parziali di tipo ellittico e parabolico:
aspetti geometrici, disuguaglianze collegate, e applicazioni}. The third author is partially supported by the GNAMPA project {\em Operatori di
Schr\"odinger con potenziali elettromagnetici singolari: stabilit\`a spettrale e stime di decadimento}. All the authors are members of the
Gruppo Nazionale per l'Analisi Matematica, la Probabilit\`a e le loro Applicazioni (GNAMPA) of the Istituto Nazionale di Alta Matematica (INdAM).

\end{document}